\documentclass[11pt]{amsart}

\usepackage{fullpage}
\usepackage{amsmath,amssymb,amsthm,amscd}

\title{{\bf Perspectives On Geometric Analysis}}
\author{Shing-Tung Yau}
\address{Department of Mathematics\\ Harvard University\\ Cambridge,
MA 02138} \email{yau@math.harvard.edu}

\thanks{To appear in {\sl Survey in Differential Geometry}, Vol. {\bf X},
2005.} \thanks{This research is supported by NSF grants
DMS-0244464, DMS-0354737 and DMS-0306600.}

\date{}

\begin{document}

\maketitle

This essay grew from a talk I gave on the occasion of the
seventieth anniversary of the Chinese Mathematical Society. I
dedicate the lecture  to the memory of my teacher S. S. Chern who
had passed away half a year before (December 2004).

During my graduate studies, I was rather free in picking research
topics. I \cite{Yau1971} worked on fundamental groups of manifolds
with non-positive curvature.  But in the second year of my
studies, I started to look into differential equations on
manifolds. However, at that time, Chern was very much interested
in the work of Bott  on holomorphic vector fields. Also he told me
that I should work on Riemann hypothesis. (Weil had told him that
it was time for the hypothesis to be settled.) While Chern did not
express his opinions about  my research on geometric analysis, he
started to appreciate it a few years later. In fact, after Chern
gave a course  on Calabi's works on affine geometry in 1972 at
Berkeley, S. Y. Cheng told me about these inspiring lectures. By
1973, Cheng and I started to work on some problems mentioned in
Chern's lectures. We did not realize that the great geometers
Pogorelov, Calabi and Nirenberg were also working on them. We were
excited that we solved some of the conjectures of Calabi on
improper affine spheres. But soon after we found out that
Pogorelov \cite{Pog1972} published his results right before us by
different arguments. Nevertheless our ideas are useful in handling
 other problems in affine geometry, and my knowledge about
Monge-Amp\`{e}re equations started to
 broaden in these years.
Chern was very pleased by my work, especially after I
\cite{Yau1978} solved the problem of Calabi on K\"{a}hler Einstein
metric in 1976. I had been at Stanford,  and Chern proposed to the
Berkeley Math Department that they hire me. I visited Berkeley in
1977 for a year and gave a course on geometric analysis with
emphasis on isometric embedding.

Chern nominated me to give a plenary talk at the International
Congress in Helsinki. The talk \cite{Yau1980} went well, but my
decision not to stay at Berkeley did not quite please him.
Nevertheless he recommended me for a position on the  faculty at
the Institute for Advanced Study. Before I accepted a faculty
position at the Institute, I organized a special year on geometry
in 1979 at the Institute  at the invitation of Borel. That was an
exciting year because most people in geometric analysis came.

In 1979, I visited China at the invitation of Professor L. K. Hua.
I gave a series of talks on the bubbling process of
Sacks-Uhlenbeck \cite{SU1981}. I suggested to the Chinese
mathematicians that they apply similar arguments for a Jordan
curve bounding two surfaces with the same constant mean curvature.
I thought it  would be a good exercise for getting into this
exciting field of geometric analysis. The problem was indeed
picked up by a group of students of Professor G. Y. Wang
\cite{JW1993}. But unfortunately it also initiated some ugly
fights during the meeting of the sixtieth anniversary of the
Chinese Mathematical Society. Professor  Wang was forced to
resign, and this event hampered  development of this beautiful
subject in China in the past ten years.

In 1980, Chern decided to develop geometric analysis on a large
scale. He initiated a series of international conferences on
differential geometry and differential equations to be held each
year in China.  For the first year, a large group of the most
distinguished mathematicians was gathered in Beijing to give
lectures (see \cite{Pro1980}). I lectured on open problems in
geometry \cite{Yau1982}. It took a much longer time than I
expected for Chinese mathematicians to pick up some of these
problems.

To his disappointment, Chern's enthusiasm about developing
differential equations and differential geometry in China did not
stimulate as much activity as he had hoped. Most Chinese
mathematicians were trained in analysis but were rather weak in
geometry. The goal of geometric analysis for understanding
geometry was not appreciated. The major research center on
differential geometry came from students of Chern, Hua and B. C.
Su. The works of J. Q. Zhong (see, e.g.,
\cite{ZhY1984,MZh1986,MZh1989}) were remarkable. Unfortunately he
died about twenty years ago. Q. K. Lu studied the Bergman metric
extensively. C. H. Gu \cite{Gu1980} studied gauge theory and
considered harmonic map where the domain is $R^{1,1}$. J. X. Hong
(see, e.g., \cite{Hong1993,HHL2003}) did some interesting work on
isometric embedding of surfaces into $\mathbb{R}^3$.  In the past
five years, the research center at the Chinese University of Hong
Kong, led by L. F. Tam and X. P. Zhu, has produced first class
work related to Hamilton's Ricci flow (see, e.g., \cite{CTZ2004,
CZ2000,CZ2005,Chen-Zhu2005,CT2005}).

In the hope that it will advance  Chern's ambition to build up
geometric analysis, I will explain my personal view to my Chinese
colleagues. I will consider this article to be successful if it
conveys to my readers the excitement of developments in
differential geometry which have been taking place during the
period when it has been my good fortune to contribute. I do not
claim this article  covers all aspects of the subject. In  fact, I
have given priority to those works closest to my personal
experience, and,  consequently, I have given insufficient space to
aspects of differential geometry in which I have not participated.
In spite of these shortcomings, I hope that my readers,
particularly those too young to know the origins of geometric
analysis, will be interested to learn how the field looks to
someone who was there. I would like to  thank comments  given by
R. Bryant,  H. D. Cao, J.  Jost,  H. Lawson, N. C. Leung, T. J.
Li, Peter Li, J. Li, K. F. Liu,  D. Phong , D. Stroock , X. W.
Wang, S. Scott, S. Wolpert and  S. W. Zhang.  I am also grateful
to J. X. Fu, especially for his help of tracking down references
for the major part of this survey. When Fu went back to China,
this task was taken up by P. Peng and X. F. Sun to whom I am
grateful also.

In this whole survey, I follow the following:

\vspace{0.5cm}
\begin{center}
{\large\bf Basic Philosophy:}
\end{center}
\vspace{0.3cm}
\begin{quote}
{\bf Functions, tensors and subvarieties governed by natural
differential equations provide deep insight into geometric
structures. Information about  these objects will give a way to
construct a geometric structure. They also provide important
information for physics, algebraic geometry and topology.
Conversely it is vital to learn ideas from these fields.

Behind such basic philosophy, there are basic invariants to
understand how space is twisted. This is provided by Chern classes
\cite{Chern1946}, which appear in every branch of mathematics and
theoretical physics. So far we barely understand the analytic
meaning of the first Chern class. It will take much more time for
geometers to understand the analytic meaning of the higher Chern
forms. The analytic expression of Chern classes by differential
forms have opened up a new horizon for global geometry. Professor
Chern's influence on mathematics is forever.}
\end{quote}
\newpage

\begin{center}
{\large \bf An old Chinese poem says:}
\end{center}

\vspace{0.8cm}

{\bf\em The reeds and rushes are abundant,\\

and the white dew has yet to dry.\\

The man whom I admire
is on the bank of the river.\\

I go against the stream in quest of him,\\

But the way is difficult and turns to the right.\\

I go down the stream in quest of him,\\

and Lo! He is on the island in the midst of the water.}

\vspace{1.5cm}

{\bf May the charm and beauty be always the guiding principle of
geometry!}

\newpage

\tableofcontents

\newpage

\section{History and contributors of the subject}

\subsection{Founding fathers of the subject}

Since the whole development of geometry depends heavily on the
past, we start out with historical developments. The following are
samples of work before 1970  which provided fruitful ideas and
methods.

\begin{itemize}
\item {\bf Fermat's principle of calculus of variation} (Shortest
path in various media).

\item {\bf Calculus (Newton and Leibnitz)}: Path of bodies
governed by law of nature.

\item {\bf Euler, Lagrange}: Foundation for the variational
principle and the study of partial differential equations.
Derivations of equations for fluids and for minimal surfaces.

\item{\bf Fourier, Hilbert}: Decomposition of functions into
eigenfunctions, spectral analysis.

\item{\bf Gauss, Riemann}: Concept of intrinsic geometry.

\item {\bf Riemann, Dirichlet, Hilbert}: Solving Dirichlet
boundary value problem for harmonic function using variational
method.

\item{\bf Maxwell}: Electromagnetism, gauge fields, unification of
forces.

\item{\bf Christoffel, Levi-Civita, Bianchi, Ricci}: Calculus on
manifolds.

\item{\bf Riemann, Poincar\'{e}, Koebe, Teichm\"uller}: Riemann surface
uniformization theory, conformal deformation.

\item{\bf Frobenius, Cartan, Poincar\'e}: Exterior differentiation
and Poincar\'e lemma.

\item{\bf Cartan}: Exterior differential system, connections on
fiber bundle.

\item{\bf Einstein, Hilbert}: Einstein equation and Hilbert
action.

\item{\bf Dirac}: Spinors, Dirac equation, quantum field theory.

\item{\bf Riemann, Hilbert, Poincar\'e, Klein, Picard, Ahlfors,
Beurling, Carlsson}: Application of complex analysis to geometry.

\item{\bf K\"{a}hler, Hodge}: K\"{a}hler metric and Hodge theory.

\item {\bf Hilbert, Cohn-Vossen, Lewy, Weyl, Hopf, Pogorelov,
Efimov, Nirenberg}: Global surface theory in three space based on
analysis.

\item{\bf Weierstrass, Riemann, Lebesgue, Courant, Douglas,
Rad\'{o}, Morrey}: Minimal surface theory.

\item{\bf Gauss, Green, Poincar\'{e}, Schauder, Morrey}: Potential
theory, regularity theory for elliptic equations.

\item{\bf Weyl, Hodge, Kodaira, de Rham,
Milgram-Rosenbloom, Atiyah-Singer}: de Rham-Hodge theory, integral
operators, heat equation, spectral theory of elliptic self-adjoint
operators.

\item{\bf Riemann, Roch, Hirzebruch, Atiyah-Singer}: Riemann-Roch
formula and index theory.

\item {\bf Pontrjagin, Chern, Allendoerfer-Weil}: Global
topological invariants defined by curvature forms.

\item{\bf Todd, Pontrjagin, Chern, Hirzebruch Grothendieck, Atiyah}:
Characteristic classes and $K$-theory in topology and algebraic
geometry.

\item {\bf Leray, Serre}: Sheaf theory.

\item{\bf Bochner-Kodaira}: Vanishing of cohomology groups based
on the curvature consideration.

\item{\bf Birkhoff, Morse, Bott, Smale}: Critical point theory,
global topology, homotopy groups of classical groups.

\item{\bf De Giorgi-Nash-Moser}: Regularity theory for the higher
dimensional elliptic equation and the parabolic equation of
divergence type.

\item{\bf Kodaira, Morrey, Grauert, Hua, H\"{o}rmander, Bergman, Kohn,
Andreotti-Vesentini}: Embedding of complex manifolds,
$\bar\partial$-Neumann problem, $L^2$ method, kernel functions.

\item{\bf Kodaira-Spencer, Newlander-Nirenberg}: Deformation of
geometric structures.

\item {\bf Federer-Fleming, Almgren, Allard, Bombieri, De Giorgi,
Giusti}: Varifolds and minimal varieties in higher dimensions.

\item{\bf Eells-Sampson, Al'ber}: Existence of  harmonic maps into
manifolds with non-positive curvature.

\item{\bf Calabi}: Affine geometry and conjectures on K\"{a}hler
Einstein metric.
\end{itemize}

\newpage

\subsection{Modern Contributors}

The major contributors can be roughly mentioned in the following
periods:

{\bf I. 1972 to 1982:} M. Atiyah, R. Bott, I. Singer, E. Calabi,
L. Nirenberg, A. Pogorelov, R. Schoen, L. Simon, K. Uhlenbeck, S.
Donaldson, R. Hamilton, C. Taubes, W. Thurston, E. Stein, C.
Fefferman, Y. T. Siu, L. Caffarelli, J. Kohn, S. Y. Cheng, M.
Kuranishi, J. Cheeger, D. Gromoll, R. Harvey, H. Lawson, M.
Gromov, T. Aubin, V. Patodi, N. Hitchin, V. Guillemin, R. Melrose,
Colin de Verdi\`{e}re, M. Taylor, R. Bryant, H. Wu, R. Greene,
Peter Li, D. Phong, S. Wolpert, J. Pitts, N. Trudinger, T.
Hildebrandt, S. Kobayashi, R. Hardt, J. Spruck, C.  Gerhardt, B.
White, R. Gulliver, F. Warner, J. Kazdan.

Highlights of the works in this period include a deep
understanding of the spectrum of elliptic operators, introduction
of self-dual connections for four manifolds, introduction of a
geometrization program for three manifolds, an understanding of
minimal surface theory, Monge-Amp\'ere equations and the
application of the theory to algebraic geometry and general
relativity.

{\bf II. 1983 to 1992:} In 1983,  Schoen and I started to give
lectures on geometric analysis at the Institute for Advanced
Study. J. Q. Zhong took notes on the majority of our lectures. The
lectures were continued in 1985 in San Diego. During the period of
1985 and 1986, K. C. Chang and W. Y. Ding came to take notes of
some part of our lectures. The book {\sl Lectures on Differential
Geometry} was published in Chinese around 1989 \cite{SY1994}. It
did have great influence for a generation of Chinese
mathematicians to become interested in this subject. At the same
time, a large group of my students made contributions to the
subject. This includes A. Treibergs, T. Parker, R. Bartnik, S.
Bando, L. Saper, M. Stern, H. D. Cao, B. Chow, W. X.  Shi, F. Y.
Zheng and G. Tian.

At the same time, J. Bismut, C. S. Lin, N. Mok, J. Q. Zhong, J.
Jost, G. Huisken, D. Jerison, P. Sarnak, K. Fukaya, T. Mabuchi, T.
Ilmanen, C. Croke, D. Stroock, Price, F. H. Lin, S. Zelditch, D.
Christodoulou, S. Klainerman, V. Moncrief, C. L. Terng, Michael
Wolf, M. Anderson, C. LeBrun, M.  Micallef, J. Moore, John Lee, A.
Chang, N. Korevaar were making contributions in various
directions. One should also mention that in this period important
 work was done  by the authors in the first group. For example,
Donaldson, Taubes \cite{Tau1982} and Uhlenbeck
\cite{Uhl1982,Uhl19822} did  spectacular work on Yang-Mills theory
of general manifolds which led Donaldson \cite{Don1983} to solve
the  outstanding question on four manifold topology. Donaldson
\cite{Don1985}, Uhlenbeck-Yau \cite{UY1986} proved the existence
of Hermitian Yang-Mills connection on stable bundles. Schoen
\cite{Sch1984} solved the Yamabe problem.

{\bf III. 1993 to now:}  Many  mathematicians joined the subject.
This includes  P. Kronheimer, B. Mrowka, J. Demailly,  T. Colding,
W. Minicozzi, T. Tao, R. Thomas, Zworski, Y. Eliashberg, Toth,
Andrews, L. F. Tam, N. C. Leung, Y. B.  Ruan, W. D. Ruan, R.
Wentworth, A. Grigor'yan, L. Saloff-Coste, J. X. Hong, X. P. Zhu,
M. T. Wang, A. K. Liu, K. F. Liu, X. F. Sun, T. J. Li, X. J. Wang,
J. Loftin, H. Bray, J. P. Wang, L. Ni, P. F. Guan, N. Kapouleas,
P. Ozsv\'{a}th, Z. Szab\'{o} and Y. I. Li. The most important
event is of course the first major breakthrough of Hamilton
\cite{Ham1997} in 1995 on the Ricci flow. I did propose to him in
1982 to use his flow to solve Thurston's conjecture. But after
this paper by Hamilton, it is finally realized that it is feasible
to solve the full geometrization program by geometric analysis. (A
key step was the estimates on parabolic equations initiated by
Li-Yau \cite{LY1986} and accomplished by Hamilton for Ricci flow
\cite{Ham1988,Ham1993}.) In 2002, Perelman \cite{Per2002,Per2003}
brought in fresh new ideas to solve important steps that are
remained in the program. Many contributors, including
Colding-Minicozzi \cite{CM2005}, Shioya-Yamaguchi
\cite{Shioya-Yamaguchi2000} and Chen-Zhu \cite{CZ2005},
\cite{Chen-Zhu2005} have helped in filling important gaps in the
arguments of Hamilton-Perelman. Cao-Zhu has just finished a long
manuscript which gives a complete detail account of the program.
The monograph will be published by International Press. In the
other direction, we see the important development of
Seiberg-Witten theory \cite{Wit1994}. Taubes
\cite{Tau19961,Tau19962,Tau19991,Tau19992} was able to prove the
remarkable theorem for counting pseudo-holomorphic curves in terms
of his invariants. Kronheimer-Mrowka \cite{KM1995} were able to
solve the Thom conjecture that holomorphic curves provide the
lowest genus surfaces in representing homology in algebraic
surfaces. (Ozsv\'{a}th-Szab\'{o} had a symplectic version
\cite{OS2000}.)

\newpage

\section{Construction of functions in geometry}\label{sec 1}   

The following is the basic principle \cite{Yau1980}:
\begin{quote}
{\noindent\bf Linear or non-linear analysis is developed to
understand the underlying geometric or combinatorial structure. In
the process, geometry will provide deeper insight of analysis. An
important guideline is that space of special  functions defined by
the structure of the space can be used to define the structure of
this space itself.}
\end{quote}

Algebraic geometers have defined the Zariski topology of an
algebraic variety using  ring of rational functions. In
differential geometry, one should  extract information about the
metric and topology of the manifolds from functions defined over
it. Naturally,  these functions should be defined either by
geometric construction or by differential equations given by the
underlying structure of the  geometry. (Integral equations have
not been used extensively as the idea of linking local geometry to
global geometry is more compatible with the ideas of  differential
equations.) A natural generalization of functions consists of the
following: differential forms, spinors, and sections of vector
bundles.

The dual concepts of differential forms or sections of vector
bundles are submanifolds or foliations. From the differential
equations that arise  from the variational principle, we have
minimal submanifolds or holomorphic cycles. Naturally the
properties of such objects or the moduli space of such objects
govern the geometry of the underlying manifold. A very good
example is Morse theory on the space of loops on a manifold (see
\cite{Mil1963}).

I shall now discuss various methods  for  constructing functions
or tensors of geometric interest.

\subsection{Polynomials from ambient space.}                        

If the manifold is isometrically embedded into Euclidean space, a
natural class of functions are the restrictions of polynomials
from Euclidean space. However, isometric embedding in general is
not rigid, and so functions constructed in such a way are usually
not too useful.

On the other hand, if a manifold is embedded into Euclidean space
in a canonical manner and the geometry of this submanifold is
defined by some group of linear transformations of the Euclidean
space, the polynomials restricted to the submanifold do play
important roles.

\subsubsection
{Linear functions being the harmonic function or eigenfunction of the
submanifold} 

For minimal submanifolds in Euclidean space, the restrictions of
linear functions are harmonic functions. Since the sum of the norm
square of the gradient of the coordinate functions is equal to
one, it is fruitful to construct classical potentials using
coordinate functions. This principle was used by Cheng-Li-Yau
\cite{CLY1984} in 1982 to give a comparison theorem for a heat
kernel for minimal submanifolds in Euclidean space, sphere and
hyperbolic space. Li-Tian \cite{LT1995} also considered a similar
estimate for complex submanifolds of $\mathbb{C}P^n$. But this
follows from \cite{CLY1984} as such submanifolds can be lifted to
a minimal submanifold in   $S^{2n+1}$.

Another very important property of a linear function is that when
it is restricted to a minimal hypersurface in a sphere $S^{n+1}$,
it is automatically an eigenfunction. When the hypersurface is
embedded, I conjectured that the first eigenfunction is linear and
the first eigenvalue of the hypersurface is equal to $n$ (see
\cite{Yau1982}). While this conjecture is not completely solved,
the work of Choi-Wang \cite{CW1983} gives strong support. They
proved that the first eigenvalue has a lower bound depending only
on $n$. Such a result was good enough for Choi-Schoen
\cite{CS1985} to prove a compactness result for embedded minimal
surfaces in $S^3$.

\subsubsection{Support functions}\label{subsec:1.1.2}                      

An important class of functions that are constructed from the
ambient space are the support functions of a hypersurface. These
are functions defined on the sphere and are related to the Gauss
map of the hypersurface. The famous Minkowski problem reduces to
solving some Monge-Amp\`{e}re equation for such support functions.
This was done by Nirenberg \cite{Nir1953}, Pogorelov
\cite{Pog1953}, Cheng-Yau \cite{CY19763}. The question of
prescribed symmetric functions of principal curvatures has been
studied by many people: Pogorelov \cite{Pog1984},
Caffarelli-Nirenberg-Spruck \cite{CNS1985}, P. F. Guan and his
coauthors (see \cite{GM2003,GG2002}), Gerhardt \cite{Ger2003},
etc. It is not clear whether one can formulate a useful Minkowski
problem for higher codimensional submanifolds.

The question of isometric embedding of surfaces into three space
can also be written in terms of the Darboux equation for the
support function. The major global result is the Weyl embedding
theorem for convex surfaces, which was proved by Pogorelov
\cite{Pogorelov1961-1,Pogorelov1961-2} and Nirenberg
\cite{Nir1953}. The rigidity part was due to Cohn-Vossen and an
important estimate was due to Weyl himself. For local isometric
embeddings, there is work by  C. S. Lin \cite{Lin1985,Lin1986},
which are followed by Han-Hong-Lin \cite{HHL2003}. The global
problem for surfaces with negative curvature was studied by  Hong
\cite{Hong1993}. In all these problems, infinitesimal rigidity
plays an important role. Unfortunately they are only well
understood for a convex hypersurface. It is intuitively clear that
generically, every closed surface is infinitesimally rigid.
However, significant works only appeared for very special
surfaces. Rado studied the set of surfaces that are obtained by
rotating a curve around an axis. The surfaces constructed depend
on the height of the curve. it turns out that such surfaces are
infinitesimally rigid except on a set of heights which form part
of a spectrum of some Sturm-Liouville operator.

\subsubsection{Gradient estimates of natural functions induced from ambient space}                                  

A priori estimates are the basic tools for nonlinear analysis. In
general the first step is to control the ellipticity of the
problem. In the case of the Minkowski problem, we need to control
the Hessian of the support function. For minimal submanifolds and
other submanifold problem, we need gradient estimates which we
shall discuss in Chapter 4. In 1974 and 1975, S. Y. Cheng and I
\cite{CY19762,CY1986} developed several gradient estimates for
linear or quadratic polynomials in order to control metrics of
submanifolds in Minkowski spacetime or affine space. This kind of
idea can be used to deal with many different metric problems in
geometry.

The first theorem concerns a spacelike hypersurface $M$ in the
Minkowski space $\mathbb{R}^{n,1}$. The following  important
question arose: Since the metric on $\mathbb{R}^{n,1}$ is
$\sum(dx^i)^2-dt^2$, the restriction of this metric on $M$ need
not be complete even though it may be true for the induced
Euclidean metric. In order to prove the equivalence of these two
concepts for hypersurfaces whose mean curvatures are controlled,
Cheng and I proved the gradient estimate of the function
$$\langle X,X\rangle=\sum_{i}(x^i)^2-t^2$$
restricted on the hypersurface.

By choosing a coordinate system, the function $\langle X,X\rangle$
can be assumed to be positive and proper on $M$. For any positive
proper function $f$ defined on $M$, if we prove the following
gradient estimate
$$\frac{\mid \bigtriangledown f\mid}{f}\leq C$$
where $C$ is independent of $f$, then we can prove the induced
metric on $M$ is complete. This is obtained by integrating the
inequality to get
$$ \mid \log f(x)-\log f(y)\mid\leq C d(x,y)$$
so that when $f(y)\rightarrow \infty$, $d(x,y)\rightarrow \infty$.
Once we knew the metric was complete, we proved the Bernstein
theorem which says that maximal spacelike hypersurface must be
linear. Such work was then generalized by Treibergs
\cite{Tre1982}, C. Gerhardt \cite{Ger1983} and R. Bartnik
\cite{Bar1984} for hypersurfaces in more general spacetime. (It is
still an important problem to understand the behavior of a maximal
spacelike hypersurface foliation for  general spacetime when we
assume the spacetime is  evolved by   Einstein equation from a
nonsingular data set.)

Another important example is the study of affine hypersurfaces
$M^n$ in an affine space $A^{n+1}$. These are the improper affine
spheres
$$\det(u_{ij})=1$$
where $u$ is a convex function or the hyperbolic affine spheres
$$\det(u_{ij})=\left(-\frac{1}{u}\right)^{n+2}$$
where $u$ is convex and zero on $\partial\Omega$ and $\Omega$ is a
convex domain. Note that these equations describe hypersurfaces
where the affine normals are either parallel or converge to a
point.

For affine geometry, there is an affine invariant metric defined
on $M$ which is
$$(\det h_{ij})^{-\frac{1}{n+2}}\sum h_{ij}dx^idx^j$$
where $h_{ij}$ is the second fundamental form of $M$. A
fundamental question is  whether this metric is complete or not.

A coordinate system in $A^{n+1}$ is chosen so that the height
function is a proper positive function defined on $M$. The
gradient estimate of the height function gives a way to prove
completeness of the affine metric. Cheng and I \cite{CY1986} did
find such an estimate which is similar to the one given above.

Once completeness of the affine metric is known, it is straight
forward to obtain important properties of the affine spheres, some
of which were conjectured by Calabi. For example we proved that an
improper affine sphere is a paraboloid and that every proper
convex cone admits a foliation of hyperbolic affine spheres. The
statement about improper affine sphere was first proved by
J\"{o}rgens \cite{Jor1954}, Calabi \cite{Cal1958} and Pogorelov
\cite{Pog1972}. Conversely, we also proved that every hyperbolic
affine sphere is asymptotic to a convex cone. (The estimate of
Cheng-Yau was reproduced again by a Chinese mathematician who
claimed to prove the result ten years afterwards.) Much more
recently, Trudinger and X. J. Wang \cite{TW2000} solved the
Bernstein problem for an affine minimal surface, thereby settling
a conjecture by Chern. They found a counterexample for
dim$\geq10$. These results are solid contributions to fourth order
elliptic equations.

The argument of using gradient estimates for some naturally
defined function was also used by me to prove that the K\"{a}hler
Einstein metric constructed by Cheng and myself is complete for
any bounded pseudo-convex domain \cite{CY1980}. (It appeared in my
paper with Mok \cite{MY1983} who proved the converse statement
which says that if the K\"{a}hler Einstein metric is complete, the
domain is pseudo-convex.)

It should be noted that in most cases, gradient estimates amount
to control of ellipticity of the nonlinear elliptic equation.

\begin{quotation}
\noindent {\bf Comment}: {To control a metric, find functions that
are capable of describing the geometry and give gradient or higher
order estimates for these functions.}
\end{quotation}

\subsection{Geometric construction of functions}                        

\subsubsection{Distance function and Busemann function.}                

When manifolds cannot be embedded into the linear space, we
construct functions adapted to the metric structure. Obviously the
distance function is the first major function to be constructed. A
very important property of the distance function is that when the
Ricci curvature of the manifold is greater than the Ricci
curvature of a model manifold which is spherical symmetric at one
point, the Laplacian of the distance function is not greater than
the Laplacian of the distance function of the model manifold in
the sense of distribution. This fact was used by Cheeger-Yau
\cite{CY1981} to give a sharp lower estimate of the heat kernel of
such manifolds. An argument of this type was also used by Perelman
in his recent work.

Gromov \cite{Gro1981} developed a remarkable Morse theory for the
distance function (a preliminary version was developed by  K.
Grove and K. Shiohama \cite{GS1977}) to compare the topology of a
geodesic ball to that of a large ball, thereby obtaining a bound
on the Betti numbers of compact manifolds with nonnegative
sectional curvature. (He can also allow the manifolds to have
negative curvature. But in this case the diameter and the lower
bound of the curvature will enter into the estimate.)

We can also take the distance function from a hypersurface and
compute the Hessian of the distance function. In general, one can
prove comparison theorems, and the principle curvatures of the
hypersurface will come into the estimates. However the upper bound
of the Laplacian of the function depends only on the Ricci
curvature of the ambient manifold and the mean curvature of the
hypersurface. This kind of calculation was used in the sixties by
Penrose and Hawking to study the focal locus of a closed surface
under the assumption that the surface is $"$trapped$"$ which means
the mean curvatures are negative in both the ingoing and the
outgoing null directions. This information allowed them to prove
the first singularity theorem in general relativity (which
demonstrates that the black hole singularity is stable under
perturbation). The distance to hypersurfaces can be used as
barrier functions to prove the existence of a minimal surface as
was shown by Meeks-Yau \cite{MY1982}, \cite{MY19822}. T. Frankel
used the idea of minimizing the distance between two submanifolds
to detect the topology of minimal surfaces. In particular, two
maximal spacelike hypersurfaces in spacetime which satisfy the
energy condition must be disjoint if they are parallel at
infinity.

Out of the distance function, we can construct the Busemann
function in the following way:
\begin{quote}
{\sl Given a geodesic ray $\gamma :[0,\infty)\rightarrow M$ so
that
$$\text{distance}(\gamma(t_1),\gamma(t_2))=t_2-t_1,$$ where
$\parallel \frac{d\gamma}{dt}\parallel=1$, one defines
 $$B_{\gamma}(x)=\lim_{\overline{t\rightarrow\infty}}(d(x,\gamma(t))-t).$$}
\end{quote}

This function generalizes the notion of a linear function. For a
hyperbolic space form, its level set  defines horospheres. For
manifolds with positive curvature, it is concave. Cohn-Vossen (for
surface) and Gromoll-Meyer \cite{GM1969} used it to prove that a
complete noncompact manifold with positive curvature is
diffeomorphic to $\mathbb{R}^{n}$.

A very important property of the Busemann function is that it is
superharmonic on complete manifolds with nonnegative Ricci
curvature in the sense of distribution. This is the key to prove
the splitting principle of Cheeger-Gromoll \cite{CG1971}. Various
versions of this splitting principle have been important for
applications to the structure of manifolds. When I \cite{Yau1978}
proved the Calabi conjecture,  the splitting principle was used by
me and others to prove the structure theorem for K\"{a}hler
manifolds with a nonnegative first Chern Class. (The argument for
the structure theorem is due to Calabi \cite{Cal1957} who knew how
to handle the first Betti number. Kobayashi \cite{Kob1961} and
Michelsohn \cite{Mic1982} wrote up the formal argument and
Beauville \cite{Bea1983} had  a survey article on this
development.)

In 1974, I was able to use the Busemann function to estimate the
volume of complete manifolds with nonnegative Ricci curvature
\cite{Yau1976}. After long discussions with me, Gromov
\cite{Gro1982} realized that my argument of Busemann function
amounts to compare volumes of geodesic balls. The comparison
theorem of Bishop-Gromov had been used extensively in metric
geometry.

If we consider $\inf_{\gamma}B_\gamma$, where $\gamma$ ranges from
all geodesic rays from a point on the manifold, we may be able to
obtain a proper exhaustion of the manifold. When $M$ is a complete
manifold with finite volume and its curvature is pinched by two
negative constants, Siu and I \cite{SiY1982} did prove that such a
function gives a concave exhaustion of the manifold. If the
manifold is also K\"{a}hler, we were able to prove that one can
compactify the manifold by adding a point to each end to form a
compact complex variety. In the other direction, Schoen-Yau
\cite{SY1982} was able to use the Busemann function to construct a
barrier for the existence of minimal surfaces to prove that any
complete three dimensional manifold with positive Ricci curvature
is diffeomorphic to Euclidean space.

The Busemann function also gives a  way to detect the angular
structure at infinity of the manifold. It can be used to construct
the Poisson kernel of hyperbolic space form. For a simply
connected complete manifold with bounded and strongly negative
curvature, it is used as  a barrier to solve the Dirichlet problem
for bounded harmonic functions, after they are  mollified at
infinity. This was achieved by Sullivan \cite{Sul1983} and
Anderson \cite{And1983}. Schoen and Anderson \cite{AS1985}
obtained the Harnack inequality for a  bounded harmonic function
and identified the Martin boundary of such manifolds. W. Ballmann
\cite{Bal1989} then studied the Dirichlet problem for manifolds of
non-positive curvature. Schoen and I \cite{SY1994} conjectured
that nontrivial bounded harmonic function exists if the manifold
has bounded geometry and a positive first eigenvalue. Many
important cases were settled in \cite{SY1994}. Lyons-Sullivan
\cite{LyoS1984} proved the existence of nontrivial bounded
harmonic functions using the non-amenability of groups acting on
the manifold.

The abundance of bounded harmonic functions on the universal cover
of a compact manifold should mean that the manifold is
``hyperbolic". Hence if the Dirichlet problem is solvable on the
universal cover, one expects the Gromov volume of the manifold to
be greater than zero.

The Martin boundary was studied by L. Ji and MacPherson (see
\cite{GJT1998,JM2002}) for the compactification of various
symmetric spaces. For product of manifolds with negative
curvature, it was determined by A. Freire \cite{Fre1991}. For rank
one complete manifolds with non-positive curvature, work has been
done by Ballman-Ledrappier \cite{BL1994} and Cao-Fan-Ledrappier
\cite{CHL2005}. It should be nice to generalize the work of L. K.
Hua on symmetric spaces with higher rank to general manifolds with
non-positive curvature. Hua found that bounded harmonic functions
satisfy extra equations (see \cite{Hua1963}).

\subsubsection{Length function defined on loop space.}                      

If we look at the space of loops in a manifold, we can take the
length of each loop and thereby define a natural function on the
space of loops. This is a function for which  Morse theory found
rich application. Bott \cite{Bott1959} made use of it to prove his
periodicity theorem. Bott \cite{Bott1954,Bott1956} and Morse also
developed formula for computing index  a geodesic. Bott showed
that the index of a closed geodesic and its linearized
Poincar\'{e} map determine the indexes of iterates of this
geodesic.  Starting from the famous works of Poincar\'{e},
Birkhoff, Morse and Ljusternik-Shnirel'man, there has been
extensive work on proving the existence of a closed geodesic using
Morse theory on the space of loops. Klingenberg and his students
developed powerful tools (see \cite{Kli1977}). Gromoll-Meyer
\cite{GM19692} did important work in which they proved the
existence of infinitely many  closed geodesics assuming the Betti
number of the free loop space of the manifold grows unboundedly.
They used the results of Bott \cite{Bott1956}, Serre and some
version of degenerate Morse theory. There was also later work by
Ballmann, Ziller, G. Thorbergsson, Hingston and Kramer (see, e.g.,
\cite{BTZ1982,Hin1984,Kra2004}), who improved the Gromoll-Meyer
theorem to give a low estimate of the growth of the number of
geometrically distinct closed geodesics of length $\leq t$. In
most cases, they grow at least as fast as the prime numbers. The
classical important question that every metric on $S^2$ supports
an infinite number of closed geodesics was also solved
affirmatively by Franks \cite{Fra1992}, Bangert \cite{Ban1993} and
Hingston \cite{Hin1993}. An important achievement was made by
Vigu\'e-Poirrier and D. Sullivan \cite{SV-P1976} who proved that
the Gromoll-Meyer condition for the existence of infinite numbers
of closed geodesics is satisfied if and only if the rational
cohomology algebra of the manifold has at least two generators.
They made use of Sullivan's theory of the rational homotopic type.
When the metric is Finsler, the most recent work of Victor Bangert
and Yiming Long \cite{BL2005} showed the existence of two closed
geodesics on the two dimensional sphere. (Katok \cite{Kat1973}
produced  an example which  shows that two is optimal.) Length
function is a natural concept in Finsler geometry. In the last
fifty years, Finsler geometry has not been popular in western
world. But under the leadership of Chern, David Bao, Z. Shen, X.
H. Mo and M. Ji did develop Finsler geometry much further (see,
e.g., \cite{BCS2000}).

A special class of manifolds,  all of whose geodesics are closed,
have occupied quite a lot of interest of distinguished geometers.
It started from the work of Zoll (1903) for surfaces where
Guillemin did important contributions. Bott \cite{Bott19542} has
determined the cohomology ring of these manifolds. The well known
Blaschke conjecture was proved by L. Green \cite{Gre1961} for two
dimension and by M. Berger and J. Kazdan (see \cite{Bes1978}) for
higher dimensional spheres. Weinstein \cite{Wei1974} and C. T.
Yang \cite{Yang1980,Yang1990,Yang1991} made important
contributions to the conjecture for other homotopic types.

\subsubsection{Displacement functions}                              

When the manifold has negative curvature, the length function of
curves is related to the displacement function defined in the
following way:

If $\gamma$ is an element of the fundamental group acting on the
universal cover of a complete manifold with non-positive
curvature, we  consider the function $d(x,\gamma(x))$: The study
of such a function gives rise to properties of compact manifolds
with non-positive  curvature. For example, in my thesis, I
generalized the Preissmann theorem to the effect that every
solvable subgroup of the fundamental group must be a finite
extension of an Abelian group which is the fundamental group of a
totally geodesic flat sub-torus \cite{Yau1971}. Gromoll-Wolf
\cite{GW1971}  and Lawson-Yau \cite{LY1972} also proved that if
the fundamental group of such a manifold has no center and splits
as a product, then the manifold splits as a metric product. Strong
rigidity result for discrete group acting on product of manifolds
irreducibly was obtained by Jost-Yau \cite{JY1997} where they
proved that these manifolds are homogeneous if the discrete group
also appears as fundamental group of compact manifolds with
nonpositive curvature.

When the manifold has bounded curvature, Margulis studied those
points where $d(x,\gamma(x))$ is  small and proved the famous
Margulis lemma which was used extensively by Gromov \cite{Gro1978}
to study the structure of manifolds with non-positive curvature.

\begin{quotation}
\noindent {\bf Comment}: {The lower bound of  sectional curvature
(or Ricci curvature) of a manifold gives upper estimate of the
Hessian (or the Laplacian) of the distance functions. Since most
functions constructed in geometry come from distance functions, we
have partial control of the Hessian of these functions. The
information provides us with basic tools to construct barrier
functions for harmonic analysis or to produce convex functions.
The Hessian of distance functions come from computations of second
variation of geodesics. If we consider the second variation of
closed geodesic loops, we get information about the  Morse index
of the loop, which enable us to link global topology to the
existence of many closed geodesics or curvatures of the manifold.

We always look for canonical objects through geometric
constructions and deform them to find their global properties.}
\end{quotation}

\subsection{Functions and tensors defined by linear differential equations}     

Direct construction of functions or tensors based on geometric
intuitions alone is not rich enough to handle the very complicated
geometric world. One should produce global geometric objects based
on global differential equations. Often the construction depends
on the maximal principle, integration by part, or the method of
contradictions, and they are not necessarily geometric intuitive.
On the other hand, basic principle of global differential
equations does fit well with modern geometry in relating local
data to global behavior. In order for the theory to be effective,
the global differential operator has to be constructed from a
geometric structure naturally.

The key to understanding any self-adjoint linear elliptic
differential operator is to understand its spectral resolution and
the detail of the structure of objects in the process of the
resolution: eigenvalues or eigenfunctions are particularly
important for their relation to geometry. Low eigenvalues and low
eigenfunctions give deep information about global geometry such as
topology or isoperimetric inequalities. High eigenvalues and high
eigenfunctions are related to local geometry such as curvature
forms or characteristic forms. Semiclassical analysis in quantum
physics give a way to relate these two ends. This results in using
 either the heat equation or the hyperbolic equation.

There are many important first order differential operators: d,
$\delta$, $\bar\partial$, Dirac operator.  All these operators
have contributed  to  a deeper  understanding of geometry. They
form systems of equations. Our understanding of them is not as
deep as the Laplacian acting on functions. The future of geometry
will rest on an understanding of global systems of equations and
their relation to global topology. The index theorem is the most
important contribution. It provides information about the kernel
(or cokernel). We still need to have a deeper understanding of the
spectrum of these operators.


\subsubsection{Laplacian}                                               

\subsubsection*{(a). Harmonic functions}
\label{subsec: 1.3.1 a - Harmonic functions}                        

The spectral resolution of the Laplacian gives rise to
eigenfunctions. Harmonic functions are therefore the simplest
functions that play important roles in geometry.

If the manifold is compact, the maximum principle shows that
harmonic functions are constant. However, when we try to
understand the singularities of compact manifolds, we may create
noncompact manifolds by scaling and blowing up processes, at which
point harmonic functions can play an important role.

The first important question about harmonic functions on a
complete manifold is  the Liouville theorem. I started my research
on analysis by understanding the right formulation of the
Liouville theorem. In 1971, I thought that it is natural to prove
that for complete manifolds with a non-negative Ricci curvature,
there is no nontrivial harmonic function \cite{Yau1975}. I also
thought that in the opposite case, when a complete manifold has
strongly negative curvature and is simply connected, one should be
able to solve Dirichlet problem for bounded harmonic functions.

The gradient estimates \cite{Yau1975} that I derived for a
positive harmonic function come from a suitable interpretation of
the Schwarz lemma in complex analysis. In fact, I generalized the
Ahlfors Schwarz lemma  before I understood how to work out the
gradient estimates for harmonic functions.  The generalized
Schwarz lemma \cite{Yau19782} says that holomorphic maps, from a
complete K\"{a}hler manifold with Ricci curvature bounded from
below to a Hermitian manifold with holomorphic bisectional
curvature bounded from above by a negative constant, are distance
decreasing with constants depending only on the bound on the
curvature. This generalization has since found many applications
such as the study of the geometry of moduli spaces by Liu-Sun-Yau
\cite{LSY2004,LSY20042}. They used it to prove the equivalence of
the Bergman metric with the K\"{a}hler-Einstein metric on the
moduli space. They also proved that these metrics are equivalent
to the Teichm\"uller metric and the McMullen metric.

The classical Liouville theorem has a natural generalization:
Polynomial growth harmonic functions are in fact polynomials.
Motivated by this fact and several complex variables, I asked
whether the space of polynomial growth harmonic functions with a
fixed growth rate is finite dimension with the  upper bound
depending only on the growth rate \cite{Yau1986}. This was proved
by Colding-Minicozzi \cite{CM1997} and generalized by Peter Li
\cite{Li1997}. (Functions can be replaced by sections of bundles).
In a beautiful series of papers (see, e.g., \cite{LW2000,
LW2002}), P. Li and J. P. Wang studied the space of harmonic
functions in relation to the geometry of manifolds. In the case
when  harmonic functions are holomorphic, they form a ring. I am
curious about the structure of this ring. In particular,  is it
finitely generated when the manifold is complete and has a
nonnegative Ricci curvature? A natural generalization of such a
question is to consider holomorphic sections of line bundles,
especially powers of canonical line bundles. This is part of
Mori's minimal model program.

\subsubsection*{(b). Eigenvalues and eigenfunctions}                      
\label{subsec: 1.3.1 b - Eigenvalues and eigenfunctions}

Eigenvalues reflect the  geometry of manifolds very precisely. For
domains, estimates of them date back to Lord Rayleigh. Hermann
Weyl \cite{Wey1912} solved a problem of Hilbert's  on the
asymptotic behavior of eigenvalues in relation to the volume of
the domain and hence initiated a new subject of spectral geometry.
P\'olya-Szeg\"o, Faber, Krahn and  Levy gave estimates of
eigenvalues of various geometric problems. On a general manifold,
Cheeger \cite{Che1970} was the first person to relate a lower
estimate of the first eigenvalue with the isoperimetric constant
(now called the Cheeger constant). One may note that many
questions on the  eigenvalue for domains are still unsolved. The
most noted one is the P\'olya conjecture which gave a sharp lower
estimate of the Dirichlet problem in terms of volume. Li-Yau
\cite{LY1983} did settle the average version of the P\'olya
conjecture.

The gradient estimate that I found for harmonic functions can be
generalized to cover eigenfunctions and Peter Li \cite{Li1979} was
the first one to apply it to finding estimates for eigenvalues for
manifolds with positive Ricci curvature. (If the Ricci curvature
has a positive lower bound, this is due to Lichnerowicz.) Li-Yau
 \cite{LY1979} then solved the well-known problem of estimating eigenvalues of
manifolds in terms of their diameter and the lower bound on their
Ricci curvature. Li-Yau conjectured the sharp constant for their
estimates, and Zhong-Yang \cite{ZhY1984} were able to prove this
conjecture by sharpening Li and Yau's arguments. A probabilistic
argument was developed by Chen and Wang \cite{CW1995} to derive
these inequalities. The precise upper bound for the eigenvalue was
first obtained by S. Y. Cheng \cite{Che1975} also in terms of
diameter and lower bound of the Ricci curvature. Cheng's theorem
gives a very good demonstration of how the analysis of functions
provides information about geometry. As a corollary of his
theorem, he proved that if a compact manifold $M^n$ has  a Ricci
curvature $\geq n-1$ and the diameter is equal to $\pi$, then the
manifold is isometric to the sphere. He used a lower estimate for
eigenvalues based on the work of  Lichnerowicz and Obata. Colding
\cite{Col1996} was able to use functions with properties close to
those  of the first eigenfunction to prove a pinching theorem
which states that: When the Ricci curvature is bounded below by
$n-1$ and the volume is close to that of the unit sphere, the
manifold is diffeomorphic to the sphere.  There is extensive work
by Colding-Cheeger \cite{CC1997,CC2000,CC20001} and Perelman (see,
e.g., \cite{BGP1992}) devoted to  the understanding of Gromov's
theory of Hausdorff convergence for manifolds. The tools they used
include the  comparison theorem, the splitting theorem of Cheeger
and Gromoll, and the ideas introduced earlier by Colding.

A very precise estimate of eigenvalues of the Laplacian has been
important in many areas of mathematics. For example, the idea of
Szeg\"{o} \cite{Sze1954}-Hersch \cite{Her1970} on the upper bound
of the first eigenvalue in terms of the area alone was generalized
by me to the higher genus in joint works with P. Yang
\cite{YY1980} and P. Li \cite{LY1982}. For genus one, this was
Berger's conjecture, as I was informed  by Cheng. After Cheng
showed me the paper of Hersch, I realized how to create trial
functions by taking the branched conformal cover of $S^2$. While
the constant in the paper of Yang-Yau \cite{YY1980}
 for torus is not the best
possible, the recent work of Jakobson, Levitin, Nadirashvili,
Nigam and Polterovich \cite{JLNNP2005} demonstrated that the
constant for a genus two surface may be the best possible and may
be achieved by Bolza's surface.

Shortly afterwards,  I  applied the argument of \cite{YY1980} to
prove that a Riemann surface defined by an arithmetic group must
have a relative high degree when it is branched over the sphere.
This observation of using Selberg's estimate coupled with Li-Yau
\cite{LY1982} was made in 1985 when I was in San Diego, where I
also used similar idea to estimate genus of mini-max surface in a
three dimensional manifolds and also to prove positivity of
Hawking mass. After I arrived in Harvard, I discussed the idea
with my colleague N. Elkies and B. Mazur.  The paper was finally
written up and published in 1995 \cite{Yau1996}. In the mean
while, ideas of using my work on eigenvalue coupled with Selberg's
work to study congruence subgroup was generalized by D. Abramovich
\cite{Abr1996} (my idea was conveyed by Elkies to him) and by P.
Zograf \cite{Zog1991} to the case where the curve has cusps. Most
recently Ian Agol \cite{Ago2005} also used similar idea to study
arithmetic Kleinian reflection groups.

 In a beautiful article, N.
Korevaar \cite{Kor1993} gave an upper bound, depending only on
genus and $n$, for the $n$-th eigenvalue $\lambda_n$ of a  Riemann
surface. His result answered a challenge of mine (see
\cite{Yau1982}) when I met him in Utah in 1989. Grigor'yan,
Netrusov and I \cite{GNY2004} were able to give a simplified proof
and apply the estimate to bound the index of minimal surfaces.
There are also works by P. Sarnak (see, e.g.,
\cite{Sar1984,IS1995}) on understanding eigenfunctions for such
Riemann surfaces. Iwaniec-Sarnack \cite{IS1995} showed that the
estimate of the maximum norm of the $n$-th eigenfunction on an
 arithmetic surface has significant interest in number theory.
Wolpert \cite{Wol1994} analyzes perturbation stability of embedded
eigenvalues and applies asymptotic perturbation theory and
harmonic map theory to show that stability is equivalent to the
non-vanishing of certain standard quantities in number theory.
There was also the work of Schoen-Wolpert-Yau \cite{SWY1979} on
the behavior of eigenvalues $\lambda_1,\cdots, \lambda_{2g-3}$ for
a compact Riemann surface of genus $g$. These are eigenvalues that
may tend to zero for metrics with curvature $-1$. However,
$\lambda_{2g-2}$, $\lambda_{2g-1}$, $\cdots$, $\lambda_{4g-1}$
always appear in $[c_g, \frac14]$ where $c_g>0$ depends only on
$g$. It will be nice to find the optimal $c_g$.

In this regard, one may mention the very deep problem of Selberg
on lower estimate of $\lambda_1$ for surfaces defined by an
arithmetic group. Selberg proved that it is grater than
$\frac{3}{16}$ and it was later improved by Luo-Rudnick-Sarnak
\cite{LRS1995}. For a higher dimensional locally symmetric space,
there is a similar question of Selberg and results similar to
Selberg's  were found by  J. S. Li \cite{Lij1991}  and
Cogdell-Li-Piatetski-Shapiro-Sarnak \cite{CLPS1991}. Many
researchers attempt to use Kazdhan's property $T$ for discrete
groups to study Selberg's problem.

There are many important properties of eigenfunctions that were
studied in the seventies. For example, Cheng \cite{Che1976} found
a beautiful estimate of multiplicities of eigenvalues of Riemann
surfaces based only on genus. The idea was used by Colin de
Verdi\`{e}re \cite{Col1986} to embedded graphes into
$\mathbb{R}^3$ when they satisfy nice combinatorial properties.
The connectivity and the topology of nodal domains are very
interesting questions. Melas \cite{Mel1992} did prove that for a
convex planar domain, the nodal line of second eigenfunctions must
intersect the boundary in exactly two points. Very little is known
about the number of nodal domains except the famous theorem of
Courant that the number of nodal domains of the $m$-th
eigenfunction is less than $m$.

There are several important questions related to the size of nodal
sets and the number of critical points of eigenfunctions. I made a
conjecture (see \cite{Yau1982}) about the area of  nodal sets, and
significant progress toward its resolution was made by
Donnelly-Fefferman \cite{DF1988}, Dong \cite{Dong1992} and F. H.
Lin \cite{Lin1991}. The number of critical points of an
eigenfunction is difficult to determine. I \cite{Yau1997} managed
to prove the existence of a critical point near the nodal set.
Jakobson and Nadirashvili \cite{JNa1999} gave a counterexample to
my conjecture  that the number of critical points of the $n$-th
eigenfunction is unbounded when $n$ tends to infinity. I believe
the conjecture is true for generic metrics and deserves to be
studied extensively. Nadirashvili and his coauthors
\cite{HHN1996,HHHN1999} were also the first to show that the
critical sets of eigenfunctions in $n$-dimensional manifold have a
finite $H^{n-2}$-Hausdorff measure. Afterwards, Han-Hardt-Lin
\cite{HHL1998} gave an explicit estimate.

When there is potential $V$, the eigenvalues of
$-\bigtriangleup+V$ are also important. When $V$ is convex, with
Singer, Wong and Stephen Yau, I applied the argument that I had
with Peter Li to estimate the gap $\lambda_2-\lambda_1$
\cite{SWYY}. When $V$ is arbitrary, I \cite{Yau2003} observed how
this gap depends on the lower eigenvalue of the Hessian of
$-\log\psi$, where $\psi$ is the ground  state. The method of
continuity was used by me in 1980 to reprove the work of
Brascamp-Lieb \cite{BL1976} on the convexity of $-\log \psi$ when
$V$ is convex (This work appeared in the appendix of \cite{SWYY}).
When $V$ is the scalar curvature, this was studied by Schoen and
myself extensively. In fact, in \cite{SY1983}, we found an upper
estimate of the first Dirichlet eigenvalue of the operator
$-\triangle+\frac12 R$ in terms of $\frac{3\pi^2}{2r^2}$ where $r$
is a certain concept of radius related to loops in a three
dimensional manifold. (If we replace  loops with higher
dimensional spheres, one can define a similar concept of radius.
It will be nice if such a concept can be related to eigenvalues of
differential forms.) This operator is naturally related to
conformal deformation, stability of minimal surfaces, etc. (The
works of D. Fischer-Colbrie and Schoen \cite{FS1980}, Micallef
\cite{Mic1984}, Schoen-Yau \cite{SY19795,SY1982} on stable minimal
surfaces all depend on an understanding of spectrum of this
operator.) The parabolic version appears in the recent work of
Perelman.

If there is a closed non-degenerate elliptic geodesic in the
manifold, Babi\v{c} \cite{Bab1968}, Guillemin and Weinstein
\cite{GW1976} found a sequence of eigenvalues of the Laplacian
which can be expressed in terms of the length, the rotation angles
and the Morse index of the geodesic.

\begin{quotation}
\noindent {\bf Comment}: {It is important to understand how
harmonic functions or eigenfunctions oscillate. Gradient estimate
is a good tool to achieve this. Gradient estimate for the log of
the eigenfunction can be used to prove the Louville theorem or
give a good estimate of eigenvalues. For higher eigenfunctions, it
is important to understand its zero set and its growth. By
controlling this information, one can estimate the dimension of
these functions. A good upper estimate for eigenvalues depends on
geometric intuition which may lead to construction of trial
functions that are more adaptive to geometry. It should be
emphasized that a clean and sharp estimate for the linear operator
is key to obtaining  good estimates for the nonlinear operator.}
\end{quotation}

\subsubsection*{(c). Heat kernel}                                     
\label{subsec: 1.3.1 c - Heat kernel}

Most of the work on the heat kernel over Euclidean space can be
generalized to those manifolds where Sobolev and  Poincar\'{e}
inequalities hold. (It should be noted that Aubin
\cite{Aub1976,Aub19762} and Talenti \cite{Tal1976} did find best
constant for various Sobolev inequalities on Euclidean space.)
These inequalities are all related to isoperimetric inequalities.
C. Croke \cite{Cro1980} was able to follow my work \cite{Yau19752}
on Poincar\'{e} inequalities to prove the Sobolev inequality
depending only on volume, diameter and the lower bound of Ricci
curvature. Arguments of John Nash were then used by Cheng-Li-Yau
\cite{CLY1981} to give estimates of the  heat kernel and its
higher derivatives. In this paper, an estimate of the injectivity
radius was derived and this estimate turns out to play a role in
Hamilton's theory of Ricci flow. A year later,
Cheeger-Gromov-Taylor \cite{CGT1982} made use of the wave kernel
to reprove this estimate. D. Stroock (see \cite{NS2005}) was able
to use his methods from Malliavin's calculus to give remarkable
estimates for the heat kernel at a pair of points where one point
is at the cut locus of another point.

The estimate of the  heat kernel was later generalized by Davies
\cite{Dav1988,Dav1989}, Saloff-Coste \cite{Sal1992} and Grigor'yan
\cite{Gri1991,Gri1995} to complete manifolds with polynomial
volume growth and volume doubling property. Since these  are
quasi-isometric invariants, their analysis  can be applied  to
graphs or discrete groups. See Grigor'yan's survey \cite{Gri1995}
and Saloff-Coste's survey \cite{Sal2005}.

On the other hand, the original gradient estimate that I derived
is a pointwise inequality that is much more adaptable to nonlinear
theory. Peter Li and I \cite{LY1986} were able to find a parabolic
version of it in 1984. We observed its significance for estimates
on the  heat kernel and its relation to the variational principle
for paths in spacetime. Coupled with the work of Cheeger-Yau
\cite{CY1981}, this gives a very precise estimate of the heat
kernel. Such ideas turn out to provide fundamental estimates which
play a crucial role in the analysis of Hamilton's Ricci flow
\cite{Ham1988,Ham1993}.

Not much is known about the heat kernel on differential forms or
differential forms with twisted coefficients. The fundamental idea
of using the  heat equation to prove the  Hodge theory came from
Milgram-Rosenbloom. The heat kernel for differential forms with
twisted coefficients does play an important role in the analytic
proof of the index theorem, as was demonstrated by
Atiyah-Bott-Patodi \cite{ABP1973}. It is the alternating sum that
exhibits cancellations and gives rise to index of elliptic
operators. When $t$ is small, the alternating sum reduces to a
calculation of curvature forms. When $t$ is large, it gives global
information on harmonic forms. Since the index of the operator is
independent of $t$, we can relate the index to characteristic
forms.

If a compact manifold is the quotient of a non-compact manifold by
a discrete group and if the heat kernel of the non-compact
manifold can be computed explicitly, it can be averaged to give
the heat kernel of the quotient manifold. Since the integral of
the later kernel on the diagonal can be computed by the spectrum
to be $\sum e^{-t\lambda_i}$, one can relate the displacement
function of the discrete group to the spectrum. This is the
Selberg trace formula relating length of closed geodesics to the
spectrum of the Laplacian.

\begin{quotation}
\noindent {\bf Comment}: {Understanding the heat kernel is almost
the same as understanding the heat equation. However, heat kernel
satisfies semi-group properties, which enables one to give a good
estimate of the maximum norm or higher order derivative norms as
long as the Sobolev inequality can be proved. It is useful to look
at the heat equation in spacetime where the Li-Yau gradient
estimate is naturally defined. The estimate provides special
pathes in spacetime for the estimate of the kernel. However, the
effects of closed geodesics have not been found in the heat
equation approach. A sharp improvement of the Li-Yau estimate may
lead to such information.}
\end{quotation}

\subsubsection*{(d). Isoperimetric inequalities}                          
\label{subsec: 1.3.1 d - Isoperimetric inequalities}

Isoperimetric inequality is a beautiful subject. It has a long
history. Besides its relation to eigenvalues, it reviews the deep
structure of manifolds. The best constant for the inequality is
important. P\'{o}lya-Szeg\"{o} \cite{PS1951}, G. Faber (1923), E.
Krahn (1925) and P. L\'{e}vy (1951)  made fundamental
contributions. Gromov generalized the idea of L\'{e}vy to obtain a
good estimate of Cheeger's constant (see \cite{Gro1999}).  C.
Croke \cite{Cro1984} and Cao-Escobar \cite{CE2003} have  worked on
domains in a simply connected manifold with non-positive
curvature. It is assumed  that the inequality holds  for any
minimal subvariety in Euclidean space. But it is not known to be
true for the best constant. Li-Schoen-Yau \cite{LSY1984} did prove
it in the case of a minimal surface with a connected boundary, and
E. Lutwak,  Deane Yang and G. Y. Zhang did some beautiful work in
the affine geometry case (see, e.g., \cite{Lut1993,LYZ2000}). In
Hamilton's proof of Ricci flow convergent to the round metric on
$S^2$, he demonstrated that the isoperimetric constant of the
metric is improving. One sees how the constant controls the
geometry of the manifold.

\begin{quotation}
\noindent {\bf Comment}: {The variational principle has been the
most important method in geometry since the Greek mathematicians.
Fixing the area of the domain and minimizing the length of the
boundary is the most classical form of isoperimetric inequality.
This principle has been generalized to much more general
situations in geometry and mathematical physics. In most cases,
one tries to prove existence of the extremal object and establish
isoperimetric inequalities by calculating corresponding quantities
for the extremal object. There is also the idea of rearrangement
or symmetrization to prove isoperimetric inequalities. In the
other direction, there is the duality principle in the alculus of
variation: instead of minimizing the length of the boundary, one
can fix it and maximize the area it encloses. The principle can be
effective in complicated variational problems.}
\end{quotation}

\subsubsection*{(e). Harmonic analysis on discrete geometry}             
\label{subsec: 1.3.1 e - Harmonic analysis on discrete geometry}

There are many other ideas in geometric analysis that can be
discretized and applied to graph theory. This is especially true
for the theory of spectrum of graphs. Some of these were carried
out by F.Chung, Grigor'yan and myself (see the reference of
Chung's survey \cite{Chu2005}). But the results in \cite{Chu2005}
are far away from being exhaustive. On the other hand, Margulis
\cite{Mar1982} and  Lubotzky-Phillips-Sarnak \cite{LPS1988} were
able to make use of discrete group and number theory to construct
expanding graphes. Methods to construct and classify these
expanding graphs are important for application in computer
science. It should be noted that Kazhdan's property (T)
\cite{Kaz1967} did play an important role in such discussions. It
is also important to see how to give a good decomposition of any
graph using the spectral method.

The most important work for the geometry of a finitely presented
group was done  by Gromov \cite{Gro1981}. He proved the
fundamental structure theorem of groups where volume grows at most
polynomially. These groups must be virtually nilpotent. Geometric
ideas were developed by Varopoulos and his coauthors
\cite{VSC1992, BCG2001} on the precise behaviors of the heat
kernel in terms of volume growth. As an application of the theory
of amenable groups, R. Brooks \cite{Bro1981} was able to prove
that if a manifold covers a compact set by a discrete group
$\Gamma$, then it has positive eigenvalue if and only if $\Gamma$
is non-amenable.

Gromov \cite{Gromov1987} also developed a rich theory of
hyperbolic groups using concepts of isoperimetric type
inequalities. It would be nice to characterize these groups that
are fundamental groups of compact manifolds with non-positive
curvature or locally symmetric spaces.

\begin{quotation}
\noindent {\bf Comment}: {The geometry of a graph or complex can
be used as a good testing ground for geometric ideas. They can be
important in understanding smooth geometric structures. Many rough
geometric concepts such as isoperimetric inequalities, can be
found on graphs and in fact they play some roles in computer
network theory. On the other hand, many natural geometric concepts
should be generalized to graphs. For example, the concept of the
fiber bundle, bundle theory over graphs and harmonic forms. It is
likely one needs to have a good way to define the concept of
equivalence between such objects. When we approximate a smooth
manifold by a graph or complex, we only care about  the limiting
object and therefore some equivalence relations should be allowed.
In the case of Cayley graph of a finitely generated group, it
depends on the choice of the generating set, and properties
independent of this generating set are preferable if we are only
interested in the group itself. In the other direction, computer
networks and other practical subjects have independent interest in
graph theory. A close collaboration between geometer and computer
scientists would be fruitful.}
\end{quotation}

\subsubsection*{(f). Harmonic analysis via hyperbolic operators}          
\label{subsec: 1.3.1 f - Harmonic analysis via hyperbolic
operators}

There are important works of Fefferman, Phong, Lieb, Duistermaat,
Guillemin, Melrose, Colin de Verdier, Taylor, Toth, Zelditch and
Sarnak on understanding the spectrum of the Laplacian from the
point of view of semi-classical analysis (see, e.g.
\cite{FP1981,DG1975,HMV2004,Sar2003}). Some of their ideas can be
traced back to  the geometric optics analysis of J. Keller. The
fundamental work of  Duistermaat and  H\"{o}rmander \cite{DH1972}
on propagation of singularities was also used extensively. There
has been a lot of progress on the very difficult question of
determining when one "Can hear the shape of a drum" by, among
others, Melrose (see \cite{Mel1996}), Guillemin \cite{Gui1996} and
Zelditch \cite{Zel2000}. (Priori to this, Guillemin and Kazhdan
\cite{Guillemin-Kazhdan1980} proved that no negatively curved
closed surface can be isospectrally deformed.) The first
counterexample for closed manifolds was given by J. Milnor
\cite{Mil1964} on a 16 dimensional torus. The idea was generalized
by Sunada \cite{Sun1985}, Gordon-Wilson \cite{GWi1997}. For
domains in Euclidean spaces, there were examples by Urakawa in
three dimensions. Two dimensional counterexamples were given by
Gordon-Webb-Wolpert \cite{GWW1992}, Wilson and  Szab\'{o}
\cite{Sza1999}. Most of the ideas for counterexamples are related
to the Selberg trace formula discussed in the section of heat
kernel. The semi-classical analysis based on the hyperbolic
operator also gives a very precise estimate or relation between
the geodesic and the spectrum. The support of the singularities of
the trace of the wave kernel $\sum e^{\sqrt{-1}t\sqrt{\lambda_i}}$
is a subset of the set of the lengths of closed geodesics. It is
difficult to achieve such results by elliptic theory. However,
most results are asymptotic in nature. It will be remarkable if
both methods can be combined.

\begin{quotation}
\noindent {\bf Comment}: {Fourier expansion has been a very
powerful tool in analysis and geometry. Practically any general
theorem in classical Fourier analysis should have a counterpart in
analysis of the spectrum of the Laplacian. The theory of geometric
optics and the propagation of a singularity gives deep
understanding of the singularity of a wave kernel. Geodesic and
closed geodesic becomes an important means to understand
eigenvalues. However, the theory has not been fruitful for the
Laplacian acting on differential forms. Should areas of minimal
submanifolds play a role? In the case of K\"ahler manifolds,
holomorphic cycles or the volume of special Lagrangian cycles
should be important, as the length of close geodesics appear in
the exponential decay term of the heat kernel. It would be useful
to sharpen the heat equation method to capture this lower order
information.}
\end{quotation}

\subsubsection*{(g). Harmonic forms}                                
\label{subsec: 1.3.1 g - Harmonic forms}

Natural generalizations of harmonic or holomorphic functions are
harmonic or holomorphic sections of bundles with connections. The
most important bundles are the exterior power of cotangent
bundles. Using the Levi-Civita connection, harmonic sections are
harmonic forms which, by the theory of de Rham and Hodge, give
canonical representation of cohomology classes. The major research
on harmonic forms comes from Bochner's vanishing theorem
\cite{Boc1948}. But our understanding is still poor except for
$1$-forms or when the curvature operator is positive, in which
case, the Bochner argument proved the manifold to be a homology
sphere. If there is any nontrivial operator which commutes with
the Laplacian, the eigenforms split accordingly. Making use of
special structures of such splitting, the Bochner method can be
more effective. For example, when the manifold is K\"{a}hler,
differential forms can be decomposed further to $(p,q)$-forms and
the  Kodaira vanishing theorem \cite{Kod1953} yields much more
powerful information, when the $(p,q)$ forms are twisted with a
line bundle or vector bundles. Similar arguments can be applied to
manifolds with a special holonomy group depending on the
representation theory of the holonomy group. When the complex
structure moves holomorphically, the subbundles of $(p,q)$ forms
in the bundle of $(p+q)$ forms do not necessarily deform
holomorphically. The concept of Hodge filtration is therefore
introduced. When we deform the complex structure around a point
where the complex structure degenerates, there is a monodromy
group acting on the Hodge filtration. The works of
Griffiths-Schmid \cite{Griffiths-Schmid1969} and Schmid's
$SL_2(\mathbb R)$ theorem \cite{Sch1973} give powerful control on
the degeneration of the Hodge structure. Deligne's theory of mixed
Hodge structure \cite{Del1971} plays a fundamental role for
studying singular algebraic varieties. The theory of variation of
Hodge structures is closely related to the study of period of the
differential forms. This theory also appears in the subject of
mirror symmetry. It is desirable to give a precise generalization
of these works to higher dimensional moduli spaces where
Kaplan-Cattani-Schmid made important contributions.

Harmonic forms give canonical representation to de Rham
cohomology. However the wedge product of harmonic forms need not
be harmonic. The obstruction comes from secondary cohomology
cooperation. K. T. Chen \cite{Che1971} studied the case of 1-forms
and Sullivan \cite{Sul1977} studied the general case and gave a
minimal model theory for a rational homotopic type of a manifold.
Using $\partial\bar{\partial}$-lemma of K\"{a}hler manifolds,
Deligne-Griffiths-Morgan-Sullivan \cite{DGMS1975} showed that the
rational homotopic type is formal for K\"{a}hler manifolds.

The importance of harmonic forms is that they give canonical
representation to the de Rham cohomology which is isomorphic to
singular cohomology over real numbers. It gives a powerful tool to
relate local geometry to global topology. In fact the vanishing
theorem of Bochner-Kodaira-Lichnerowicz allows one to deduce from
sign of curvature to vanishing of cohomology. This has been one of
the most powerful tools in geometry in the past fifty years.

The idea of harmonic forms came from fluid dynamics and Maxwell
equations. The non-Abelian version is the Yang-Mills theory. Most
of the works on Yang-Mills theory have been focused on these gauge
fields where the absolute minimum is achieved by some
(topological) characteristic number. (These are called BPS state
in physics literature.) When the dimension of the manifold is
four, the star operator maps two form to two form and it makes
sense to require the curvature form to be self-dual or
anti-self-dual depending whether the curvature form is invariant
or anti-invariant under the star operator. These curvature forms
can be interpreted as non-Abelian harmonic forms. The remarkable
fact is that when the metric is K\"ahler, the anti-self-dual
connections give rise to holomorphic bundles. The moduli space of
such bundles can often be computed using tools from algebraic
geometry.

If we take the product space $M\times \mathcal M$ where $M$ is the
four dimensional manifold and $\mathcal M$ is the moduli space of
anti-self-dual connections, there is a universal bundle $V$ over
$M\times\mathcal M$. By studying the slant product and the Chern
classes of $V$, we can construct polynomials on the cohomology of
$M$ that are invariants of the differentiable structure of $M$.
These are Donaldson polynomials (see \cite{DKr1990}). In general
$\mathcal M$ is not compact and Donaldson has to construct cycles
in $\mathcal M$ for such operations. Donaldson invariants are
believed to be equivalent to Seiberg-Witten invariants, where the
vanishing theorem can apply and powerful geometric consequences
can be found. Kronheimer and Mrowka \cite{KM1995} built an
important concept of  simple type for Donaldson invariants. It is
believed that Donaldson invariants of algebraic surfaces of
general type are of simple type.

If the manifold is symplectic, we can look at the moduli space of
pseudo-holomorphic curves. (These are $J$-invariant maps from
Riemann surfaces to the manifold. $J$ is an almost complex
structure that is tame to the symplectic form.) Symplectic
invariants can be created and they are called Gromov-Witten
invariants. Y. Ruan \cite{Ruan1996} has observed that they need
not be diffeomorphic invariants. It may still be interesting to
know whether Gromov-Witten invariants are invariants of
differentiable structures for Calabi-Yau manifolds.

De Rham cohomology can only capture the non-torsion part of the
singular cohomology.  Weil \cite{Wei1952} and Allendoerfer-Eells
\cite{AE1958} attempted to use differential forms with poles to
compute cohomology with integer coefficients. Perhaps one should
study Chern forms of a complex bundle with a connection that
satisfies the Yang-Mills equation and whose curvature is square
integrable. The singular set of the connection may be  allowed  to
be minimal submanifolds. The moduli space of such objects may give
information about  integral cohomology. It should be noted that
Cheeger-Simons \cite{CSi1983} did develop a rich theory of
differential character with values in $\mathbb R / \mathbb Z$. It
depends on the connections of the bundle. Witten managed to
integrate the Chern-Simons forms \cite{CS1974} on the space of
connections to obtain the knot invariants of Jones \cite{Jon1987}.

When we look for different operators acting on different forms, we
may have to look into different kinds of harmonic forms. For
example, if we are interesting in $\partial\bar\partial$
cohomology, we may look for the operator $\left (
\partial\bar\partial \right )^\ast \partial\bar\partial
+\partial\partial^\ast + \bar\partial\bar\partial^\ast$. It would
be interesting to see how super-symmetry may generalize the
concept of harmonic forms.

\begin{quotation}
\noindent {\bf Comment}: {The theory of harmonic form is
tremendously powerful because it provides a natural link between
global topology, analysis, geometry, algebraic geometry and
arithmetic geometry. However, our analytic understanding of high
degree forms is poor. For one forms, we can integrate along paths.
For two forms, we can take an interior product with a vector field
to create a moment map. For closed $(1,1)$-forms in a K\"ahler
manifold, we can express them locally as $\partial\bar\partial f$.
However, we do not have good ways to reduce a high degree form to
functions which are easier to understand. Good estimates of higher
degree forms will be very important.}
\end{quotation}

\subsubsection{$\bar{\partial}$-operator}                     
\label{subsec: 1.3.2}

Construction of holomorphic functions or holomorphic sections of
vector bundles and holomorphic curves are keys to understanding
complex manifolds.

In order to demonstrate the idea  behind the philosophy of
determining the structure of manifolds by function theory, I was
motivated to generalize the uniformization theory of a Riemann
surface to higher dimensions when I was a graduate student. During
this period, I was influenced by the works of Greene-Wu
\cite{GW1977} in formulating these conjectures. Greene and Wu were
interested in knowing whether the manifolds are Stein or not.

When the complete K\"{a}hler manifold is compact with positive
bisectional curvature, this is the Frankel conjecture, as was
proved independently by Mori \cite{Mori1979} and Siu-Yau
\cite{SiY1980}. Both arguments depend on the construction of
rational curves of low degree. Mori's argument is stronger, and it
will be good to capture his result by the analytic method. When
the manifold has nonnegative bisectional curvature and positive
Ricci curvature, Mok-Zhong \cite{MZh1986} and Mok \cite{Mok1988},
using ideas of Bando \cite{Ban1984} in his thesis on Hamilton's
Ricci flow, proved that the manifold is Hermitian symmetric.

When the complete K\"{a}hler manifold is noncompact with positive
bisectional curvature, I conjectured that it must be biholomorphic
to $\mathbb{C}^n$ (see \cite{Yau1982}). Siu-Yau \cite{SiY1977}
made the first attempt to prove such a conjecture by using the
$L^2$-method of H\"{o}rmander \cite{Hor1990} to construct
holomorphic functions with slow growth. (Note that H\"{o}rmander's
method goes back to Kodaira, which was also generalized by
Calabi-Vesentini \cite{CV1959}.) Singular weight functions were
used in this paper and later much more refined arguments were
developed by Nadel \cite{Nad1990} and Siu \cite{Siu2001} using
what is called the multiplier ideal sheaf method. Siu found
important applications of this method in algebraic geometry and
also related the idea to the powerful work of J. Kohn on weakly
pseudo-convex domain.

This work of Siu-Yau was followed by Siu-Mok-Yau \cite{MSY1981}
and Mok \cite{Mok1983,Mok1984} under assumptions about  the decay
of curvature and volume growth. Shi \cite{Shi1990,Shi1997,Shi1998}
introduced Hamilton's Ricci flow to study my conjecture, and his
work is fundamental. This was followed by beautiful works of Cao
\cite{Cao1992,Cao1997}, Chen-Zhu \cite{CZ2003,CZ2002} and
Chen-Tang-Zhu \cite{CTZ2004}. Assuming the manifold has maximal
Euclidean volume growth and bounded curvature, Chen-Tang-Zhu
\cite{CTZ2004} (for complex dimension two) and then Ni
\cite{Ni2005} (for all higher dimension) were able to prove the
manifold can be compactified as a complex variety. Last year,
Albert Chau and Tam \cite{CT2005} were finally able to settle the
conjecture assuming maximal Euclidean volume growth and bounded
curvature. An important lemma of L. Ni \cite{Ni2005} was used,
where a conjecture of mine (see \cite{Yau1991} or the introduction
of \cite{Ni2005}) was proved. The conjecture says that maximal
volume growth implies scalar curvature decays quadratically  in
the average sense.

While we see  great accomplishments for K\"ahler manifolds with
positive curvature, very little is known for K\"{a}hler manifolds,
which are complete simply connected with strongly negative
curvature. It is conjectured to be a bounded domain in
$\mathbb{C}^n$. (Some people told me that Kodaira considered a
similar problem. But I cannot find the appropriate reference.) The
major problem is to construct bounded holomorphic functions.

The difficulty of construction of bounded holomorphic functions is
that the basic principle of the $L^2$-method of H\"{o}rmander
comes from Kodaira's vanishing theorem. It is difficult to obtain
elegant results by going from weighted $L^2$ space to bounded
functions. In this connection, I was able to show that non-trivial
bounded holomorphic functions do not exist on a complete manifold
with non-negative Ricci curvature \cite{Yau19782}.

If the manifold is the universal cover of a compact K\"{a}hler
manifold $M$ which has a homotopically nontrivial map to a compact
Riemann surface with genus $>1$, then one can construct a bounded
holomorphic function, using arguments of Jost-Yau \cite{JY1992}.
In particular, if $M$ has a  map to a product of Riemann surfaces
with genus $>1$ with nontrivial topological degree, the universal
cover should have a good chance to be a bounded domain.

Of course, this kind of construction is based on the fact that
holomorphic functions are harmonic. Certain rigidity based on
curvature forced the converse to be true. For functions, the
target space has no topology and rigidity is not expected. Bounded
holomorphic functions can not be constructed by solving the
Dirichlet problem unless some boundary condition is assumed. This
would make good sense if the boundary has a nice CR structure.
Indeed, for odd dimensional real submanifold in $\mathbb{C}^n$
which has maximal complex linear subspace on each tangent plane,
Harvey-Lawson \cite{HL1975,HL1977} proved the remarkable theorem
that they bound complex submanifolds. Unfortunately the boundary
of a complete simply connected manifold with bounded negative
curvature does not have a smooth boundary. It will be nice to
define a CR structure on such a singular boundary. One may mention
the remarkable work of Kuranishi \cite{Kur19821,Kur19822,Kur19823}
on embedding of an abstract CR structure.

Historically a  motivation for the development of the
$\bar\partial$ operator came from the Levi problem, which was
solved by Morrey, Grauert and greatly improved by Kohn and
H\"ormander. Their methods are powerful in studying pseudoconvex
manifolds.

In this regard, one may mention the conjecture of Shafarevich that
the universal cover of an algebraic manifold is pseudoconvex. Many
years ago, I conjectured that if the second homotopy group of the
manifold is trivial, its universal cover can be embedded into a
domain of some algebraic manifold where the covering
transformations act on the domain by birational transformations.
One may also mention the work of S. Frankel \cite{Fra1989} on
proving that an algebraic manifold is Hermitian symmetric if the
universal cover is a convex domain in complex Euclidean space.

\begin{quotation}
\noindent {\bf Comment}:  {The $\bar\partial$ operator is the
fundamental operator in complex geometry. Classically it was used
to solve the uniformization theorem, the Levi problem and the
Corona problems. We have seen much progress on the higher
dimensional generalizations of the first two problems. However,
due to poor understanding of the construction of bounded
holomorphic functions, we are far away from understanding the
Corona problem in a higher dimension and many related geometric
questions.}
\end{quotation}

\subsubsection{Dirac operator}                                     
\label{subsec: 1.3.3 - dirac operator}

A very important bundle is the bundle of spinors. The Dirac
operator acting on spinors is the most mysterious but major
geometric operator. Atiyah-Singer were the first mathematicians to
study it in geometry and by thoroughly understanding the Dirac
operator, they were able to prove their celebrated index theorem
\cite{AS1963}. On a K\"{a}hler manifold, the Dirac operator can be
considered as a $\bar{\partial}+\bar{\partial}^*$ operator acting
on differential forms with coefficients on the square root of the
canonical line bundle. Atiyah-Singer's original proof can be
traced back to the celebrated Riemann-Roch-Hirzebruch formula and
the Hirzebruch index formula. The formulas of Gauss-Bonnet-Chern
and Atiyah-Singer-Hirzebruch should certainly be considered as the
most fundamental identities in geometry. The vanishing theorem of
Lichnerowicz \cite{Lic1963} on harmonic spinors over spin
manifolds with positive scalar curvature gives strong information.
Through the Atiyah-Singer index theorem, it gives the  vanishing
theorem for the $\hat{A}$-genus and the $\alpha$ invariants for
spin manifolds with positive scalar curvature. The method was
later sharpened by Hitchin \cite{Hit1974} to prove that every
Einstein metric over $K3-$surfaces must be K\"{a}hler and Ricci
flat. An effective use of Lichnerowicz formula for a
$\text{spin}_{\mathbb{C}}$ structure for a four dimensional
manifold is important for Seiberg-Witten theory, which couples the
Dirac operator with a complex line bundle. Lawson-Yau
\cite{LY1974} were able to use Lichnerowicz's work coupled with
Hitchin's work to prove  a large class of smooth manifolds have no
smooth non-Abelian group action and, by using modular forms, K. F.
Liu proved a loop space analogue of the Lawson-Yau's theorem for
the vanishing of the Witten genus in \cite{Liu1995}.

On the basis of the  surgery result of Schoen-Yau
\cite{SY19795,SY19791} and Gromov-Lawson \cite{GL19801,GL19802},
one expects that a suitable converse to Lichnerowicz's theorem
exists. The chief result is that surgery on spheres with
codimension $\geq 3$ preserves a class of metrics with positive
scalar curvature. Once geometric surgery is proved, standard works
on cobordism theory allow one to deduce  existence results for
simply connected manifolds with positive scalar curvature. The
best work in this direction is due to Stolz \cite{Sto1992} who
gave a complete answer in the case of  simply connected manifolds
with dimension greater then $4$. I also suggested the possibility
of performing  surgery on an asymptotic hyperbolic manifold with
conformal boundary whose scalar curvature is positive. This is
related to the recent work of Witten-Yau \cite{WY2000} on the
connectedness of the conformal boundary.

The study of metrics with positive scalar curvature is the first
important step in understanding  the positive mass conjecture in
general relativity. Schoen-Yau \cite{SY19792,SY1981} gave the
first proof using ideas of minimal surfaces. Three years later,
Witten \cite{Wit1981} gave a proof using harmonic spinors. Both
approaches have been fundamental to questions related to mass and
other conserved quantities in general relativity. In the other
direction, Schoen-Yau \cite{SY19791} generalized their argument in
1979 to find topological obstructions for higher dimensional
manifolds with positive scalar curvature. Subsequently
Gromov-Lawson \cite{GL19801,GL19802} observed that the
Lichnerowicz theorem can be coupled with a fundamental group and
give topological obstructions for a  metric with positive scalar
curvature.  This work was related to the Novikov conjecture where
many authors, including Lusztig \cite{Lus1972}, Rosenberg
\cite{Ros1983}, Weinberger \cite{Wei1990} and G. L. Yu
\cite{Yu1998} made contributions.

Besides its importance in demonstrating the stability of Minkowski
spacetime, the positive mass conjecture was used by Schoen
\cite{Sch1984} in a remarkable manner to finish the proof of the
Yamabe problem where Trudinger \cite{Tru1968} and Aubin
\cite{Aub19760} made substantial contributions.

\begin{quotation}
\noindent {\bf Comment}: {The Dirac operator is perhaps one of the
most mysterious operators in geometry. When it is twisted with
other bundles, it gives the symbol of all first order elliptic
operators. When it couples with a complex line bundle it gives the
Seiberg-Witten theory which provides powerful information for four
manifolds. On the other hand, there were two different methods to
study metrics with positive scalar curvature. It should be
fruitful to  combine both methods: the method of Dirac operator
and the method of minimal submanifolds.}
\end{quotation}

\subsubsection{First order operator twisted by vector fields or endomorphisms of
bundles}

Given a vector field $X$ on a manifold, we can consider the
complex of differential forms $\omega$ so that $L_X\omega=0$. On
such complex, $d+\iota_X$ defines a differential and the resulting
cohomology is called equivariant cohomology.

During the seventies, Bott \cite{Bott1967} and Atiyah-Bott
\cite{AB1984} developed the localization formula for equivariant
cohomology. Both the concepts of a moment map and equivariant
cohomology have become very important tools for computations of
various geometric quantities, especially Chern numbers of natural
bundles. The famous work of Atiyah, Guillemin-Sternberg on the
convexity of the image of the moment map gives a strong
application of equivariant cohomology to toric geometry. The
formula of Duistermaat-Heckman \cite{DH1982} played an important
role in motivation for evaluation of path integrals.  These works
have been used by Jeffrey and Kirwan \cite{JK1997} and by K. F.
Liu and his coauthors on several topics: the mirror principle
(Lian-Liu-Yau \cite{LLY1997,LLY1999,LLY19992,LLY2000}),
topological vertex (Li-Liu-Liu-Zhou \cite{LLLZ}), etc. The idea of
applying localization to enumerative geometry was initiated  by
Kontsevich \cite{Kon1994} and later by Givental \cite{Giv1996} and
Lian-Liu-Yau \cite{LLY1997} independently. (Lian-Liu-Yau
\cite{LLY1997} formulated a functorial localization formula which
has been fundamental for various calculations in mirror geometry.)
These works solve the identities conjectured by Candelas et al
\cite{COGP1991} based on mirror symmetry, and provide  good
examples of the ways in which conformal field theory can be a
source of inspiration when looking at classical problems in
mathematics.

If we twist the $\bar\partial$ operator with an endomorphism
valued holomorphic one form $s$ so that $s\circ s=0$, it gives
rise to a complex $\bar\partial+s$. This was the Higgs theory
initiated by Hitchin \cite{Hit19742} and studied extensively by
Simpson \cite{Sim1988}. There is extensive work of Zuo Kang and
Jost-Zuo (see \cite{Zuo1999}) on Higgs theory and representation
of fundamental groups of algebraic manifold.

In string theory, there is a three form $H$ and the cohomology of
$d^c+H$ has not been well understood. It would be interesting to
develop a deeper understanding of such twisted cohomology and its
localization.

\begin{quotation}
\noindent {\bf Comment}: {The idea of deforming a de Rham operator
by twisting with some other zero order operators has given
powerful information to geometry. Witten's idea of the analytic
proof of Morse theory is an example. Equivariant cohomology is
another example. We expect to see more works in such directions. }
\end{quotation}

\subsubsection{Spectrum and global geometry}                   
\label{subsec: 1.3.4 - Spectrum and global geometry}

Weyl made a famous address in the early fifties. The title of his
talk was  {\sl The Eigenvalue Problem Old and New}. He was excited
by the work of Minakshisundaram and Pleijel which  asserts that
the zeta function $\zeta(s)=\sum_\lambda\lambda^{-s}$, where
$\lambda$ are eigenvalues of the Laplacian, not only makes  sense
for Re's large, but also has meromorphic extension to the whole
complex $s$-plane, the position of whose poles could be described
explicitly. In particular, it is analytic near $s=0$. Formally
$\frac{d\zeta(s)}{ds}\mid_{s=0}$ can be viewed as $-\log \det
(\Delta)$. This gives a definition of determinant of Laplacian
which entered into the fundamental work of Ray-Singer relating
Reidemeister's combinational invariant of a manifold with analytic
torsion defined by determinants of the Laplacians acting on
differential forms of various degrees. Other application of zeta
function expressed in terms of kernel is the calculation of the
asymptotic growth of eigenvalues in terms of volume of the
manifold. Tauberian type theorem is needed.

This initiated the subject of finding formula to relate spectrum
of manifolds with their global geometry. Atiyah and Singer
\cite{AS1963} were the most important contributors to this
beautiful subject. Atiyah-Bott-Patodi \cite{ABP1973} applied the
heat kernel expansion to a proof of the local index theorem.
Atiyah-Patodi-Singer \cite{APS19751,APS19752,APS1976} initiated
the study of spectrum flow and gave important global spectral
invariants on odd dimensional manifolds. These global invariants
become boundary terms for the $L^2$-index theorem developed by
Atiyah-Donnelly-Singer \cite{ADS1983} and Mark Stern
\cite{Ste1989}. (A method of Callias \cite{Cal1978} has been used
for such calculations.) Witten \cite{Wit1983,Wit1984} has
introduced supersymmetry and analytic deformation of the  de Rham
complex to Morse theory, and thereby revealed a new aspect of the
connection between global geometry and theoretical physics.
Witten's work has been generalized by Demailly \cite{Dem1991} and
Bismut-Zhang \cite{BZ1992,BZ1994} to study the holomorphic Morse
inequality and analytic torsion. Novikov \cite{Nov1981} also
studied Morse theory for one forms. Witten's work on Morse theory
inspired the work of Floer (see, e.g.,
\cite{Flo19881,Flo19882,Flo1989}) who used his ideas in Floer
cohomology to prove Arnold's conjecture in case where the manifold
has vanishing higher homotopic group. Floer's theory is related to
knot theory (through Chern-Simon's theory \cite{CS1974}) on three
manifolds. Atiyah, Donaldson, Taubes , Dan Freed, P. Braam, and
others (see, e.g., \cite{Ati1988,Tau1990,BD1995,Fre1995}) all
contributed to this subject. Fukaya-Ono \cite{FO1999}, Oh
\cite{Oh1995}, Kontsevich \cite{Kon19942}, Hofer-Wysocki-Zehnder
\cite{HWZ1998}, G. Liu-Tian \cite{gLT1998}, all studied such a
theory in symplectic geometry. Some part of Arnold's conjecture on
fixed points of groups acting on symplectic manifolds was claimed
to be proven. But a completely satisfactory proof has not been
forthcoming.

One should also mention here the very important work of Cheeger
\cite{Che1977} and M\"{u}ller \cite{Mul1978} in which they verify
the conjecture of Ray-Singer equating analytic torsion with the
combinational torsion of the manifold. The fundamental idea of
Ray-Singer \cite{RS1973} on holomorphic torsion is still being
vigorously developed. It appeared in the beautiful work of Vafa et
al \cite{BCOV1994}. Many more works on analytic torsion were
advanced by Quillen, Todorov, Kontsevich, Borcherds , Bismut,
Lott, Zhang and Z. Q. Lu (see \cite{BG2004} and it's reference,
\cite{JT1998}, \cite{Bor1995,Bor1996}). The local version of the
index theorem by Atiyah-Bott-Patodi \cite{ABP1973} was later
extended in an sophisticated way by Bismut \cite{Bis1985} to an
index theorem for a family of elliptic operators.(The local index
argument dates back to the foundational work of McKean-Singer
\cite{MS1967} where methods were developed to calculate
coefficients of heat kernel expansion.) The study of elliptic
genus by Witten \cite{Wit1987}, Bott-Taubes \cite{BT1989}, Taubes
\cite{Tau1989}, K. F. Liu \cite{Liu1996} and M. Hopkins
\cite{Hop1995} has built a bridge between topology and modular
form.

\begin{quotation}
\noindent {\bf Comment}: {The subject of relating the spectrum to
global topology is extremely rich. It is likely that we have only
touched part of this rich subject. The deformation of spectrum
associated with the  deformation of geometric structure is always
a fascinating subject. Global invariants are created by spectral
flows. Determinants of elliptic operators are introduced to
understand measures of infinite dimensional space. Geometric
invariants that are created by asymptotic expansion of heat or
wave kernels are in general not well understood. It will be a long
time before we have a much better understanding of the global
behavior of spectrum.}
\end{quotation}

\newpage

\section{Mappings between manifolds and rigidity of geometric structures}
\label{sec: 2}

There is a need to exhibit a geometric structure in a simpler
space: hence we embed algebraic manifolds into complex projective
space, we isometrically embed a Riemannian manifold into Euclidean
space and we classify structures such as bundles by studying maps
into Grassmannian.

We are also interested in probing the structure of a manifold by
mapping Riemann surfaces inside the manifold, an important example
being holomorphic curves in algebraic manifolds. Of course, we are
also interested in maps that can be used to compare the geometric
structures of different manifolds.

\subsection{Embedding}

\subsubsection{Embedding theorems}                                 
\label{subsec: 1.3.5 - Embedding theorems}

Holomorphic sections of holomorphic line bundles have always been
important in algebraic geometry. The Riemann-Roch formula coupled
with vanishing theorems gave very powerful existence results for
sections of line bundles. The Kodaira embedding theorem
\cite{Kod1954}  which said that every Hodge manifold is projective
has initiated the theory of holomorphic embedding of K\"{a}hler
manifolds. For example, Hirzebruch-Kodaira \cite{HK1957} proved
that every odd (complex) dimensional K\"{a}hler manifold
diffeomorphic to projective space is biholomorphic to projective
space. (I proved the same statement  for even dimensional
K\"{a}hler manifolds based on K\"{a}hler Einstein metric.)

Given an orthonormal basis of holomorphic sections of a very ample
line bundle, we can embed the manifold into projective space. The
induced metric is the Bergman metric associated with the line
bundle. Note that the original definition of the Bergman metric
used the canonical line bundle and $L^2$-holomorphic sections.

In the process of understanding the relation between stability of
a manifold and the existence of the K\"{a}hler Einstein metric, I
\cite{Yau1986} proposed that every Hodge metric can be
approximated by the Bergman metric as long as we allow the power
of the line bundle to be  large. Following the ideas of the paper
of Siu-Yau \cite{SiY1977}, Tian \cite{Tian1990} proved the $C^2$
convergence in his thesis under my guidance. My other student W.
D. Ruan \cite{Rua1998} then proved $C^\infty$ convergence in his
thesis. This work was followed by Lu \cite{Lu2000}, Zelditch
\cite{Zel1998} and Catlin \cite{Cat1997} who observed that the
asymptotic expansion of the kernel function follows  from some
rather standard expressions of the Szeg\"{o} kernel, going back to
Fefferman \cite{Fef1974} and Boutet de Monvel-Sj\"{o}strand
\cite{BS1976} on the circle bundle associated with the holomorphic
line bundle over the K\"{a}hler manifold. Recently, Dai, Liu, Ma
and Marinescu \cite{DLM2004} \cite{MM2004} obtained the asymptotic
expansion of the kernel function by using the heat kernel method,
and gave a general way to compute the coefficients, thus also
extended it to symplectic and orbifold cases.

Kodaira's proof of embedding Hodge manifolds by the sufficiently
high power of a positive line bundle is not effective. Matsusaka
\cite{Mat1972} and later Koll\'{a}r \cite{Kol1993}, Siu
\cite{Siu1996} were able to provide effective estimate of the
power. Demailly \cite{Dem1993,Dem1996} and Siu
\cite{Siu1996,Siu2001} made a remarkable contribution toward the
solution of the famous Fujita conjecture \cite{Fuj1987} (see also
Ein and Lazarsfeld \cite{EL1993}). Siu's powerful method also
leads to a proof of the deformation invariance of plurigenera of
algebraic manifolds \cite{Siu1998}. It should be noticed that the
extension theorem of Ohsawa-Takegoshi played an important role in
this last work of Siu.

\begin{quotation}
\noindent {\bf Comment}: {The idea of embedding a geometric
structure is clearly important as once they are put in the same
space, we can compare them and study the moduli space of the
geometric structure much better. For example, one can define Chow
coordinate of a projective manifold and we can study various
concepts of geometric stability of these structures. At this
moment, there is no natural universal space of K\"ahler manifolds
or complex manifolds as we may not have a positive holomorphic
line bundle over such manifolds to embed into complex projective
space. In a similar vein, it will be nice to find a universal
space for symplectic manifolds.}
\end{quotation}

\subsubsection{Compactification}                       
\label{subsec: 1.3.6 - compactification}

Problem of compactification of the manifold dates back to Siegel,
Satake, Baily-Borel \cite{BB1966} and Borel-Serre \cite{BoS1973}.
They are important for representation theory, for algebraic
geometry and for number theory.

For geometry of non-compact manifolds, we like to control behavior
of differentiable forms at infinity. A good exhaustion function is
needed.

Construction of a proper exhaustion function with a bounded
Hessian on a complete manifold with a bounded curvature was
achieved by Schoen-Yau \cite{SY1994} in 1983 in our lectures in
Princeton. Based on this exhaustion function, M. Dafermos
\cite{Daf1997} was able to reprove a theorem of Cheeger-Gromov
\cite{CG1991} that such manifolds admit an exhaustion by compact
hypersurfaces with bounded second fundamental form. Such
exhaustions are useful to understand characteristic forms on
noncompact manifolds as the boundary term can be controlled by the
second fundamental form of the hypersurfaces.

 After my work with Siu \cite{SiY1982} on compactification
of a strongly, negatively curved K\"{a}hler manifold with finite
volume, I proposed that every complete K\"{a}hler manifold with
bounded curvature, finite volume and finite topology should be
compactifiable to be a compact complex variety. I suggested this
problem to Mok and Zhong in 1982 who did significant work
\cite{MZh1989} in this direction. (The compactification by
Mok-Zhong is not canonical and it is desirable to find an
algebraic geometric analogue of Borel-Baily compactification
\cite{BB1966} so that we can study the $L ^2$-cohomology in terms
of the intersection cohomology of the compactification.) Recall
that  the important conjecture of Zucker on identifying
$L^2$-cohomology with the intersection cohomology of the
Borel-Baily compactification was settled by Saper-Stern
\cite{SS1990} and Looijenga \cite{Loo1988}. (Intersection
cohomology was introduced by Goresky-MacPherson
\cite{GM1980,GM1983}. It is a topological concept and hence the
Zucker conjecture gives a topological meaning of the
$L^2$-cohomology.) It would be nice to find compactification for
algebraic varieties so that suitable form of intersection
cohomology can be used to understand $L^2$ cohomology.
Goresky-Harder-MacPherson \cite{GHM1994} and Saper \cite{Sap2005}
have contributed a lot toward this kind of question. For moduli
space of bundles, or polarized projective structures,
compactification means degeneration of these structures in a
suitable canonical manner. For algebraic curves, there is
Deligne-Mumford compactification \cite{DMu1969} which has played a
fundamental role in understanding algebraic curves. Geometric
invariant theory (see \cite{MFK1994}) gives a powerful method to
introduce the concept of stable structures. Semi-stable structures
can give points at infinity. The compactification based on the
geometric invariant theory for moduli space of surfaces of the
general type was done by Gieseker \cite{Gie19772}. For a higher
dimension, this was done by Viehweg \cite{Viehweg1995}. Detailed
analysis of the divisors at infinity is still missing.

\begin{quotation}
\noindent {\bf Comment}: {Compactification of a manifold is very
much related to the embedding problem. One needs to construct
functions or sections of bundles near infinity. For the moduli
space of geometric structures, it amounts to degenerately the
structures canonically. It will be important to study the
degeneration of Hermitian Yang-Mills connections and K\"{a}hler
Einstein metrics.}
\end{quotation}

\subsubsection{Isometric embedding}

Given a metric tensor on a manifold, the problem of isometric
embedding is equivalent to find enough functions $f_1, \cdots,
f_N$ so that the metric can be written as $\sum(df_i)^2$. Much
work was accomplished for two dimensional surfaces as was
mentioned in section \ref{subsec:1.1.2}. Isometric embedding for
the  general dimension was solved in the famous work of J. Nash
\cite{Nash1954, Nash1956}. Nash used his famous implicit function
theorem which depends on various smoothing operators to gain
derivatives. In a remarkable work, G\"unther \cite{Gunther1989}
was able to avoid the Nash procedure. He used only the standard
H\"older regularity estimate for the Laplacian to reproduce the
Nash isometric embedding with the same regularity result. In his
book \cite{Gromov1986}, Gromov was able to lower the codimension
of the work of Nash. He called his method the $h$-principle.

When the dimension of the manifold is $n$, the expected dimension
of the Euclidean space for the manifold to be isometrically
embedded is $\frac{n(n+1)}{2}$. It is important to understand
manifolds isometrically embedded into Euclidean space with this
optimal dimension. Only in such a dimension does  it make sense to
talk about rigidity questions. It remains a major open problem
whether one can find a nontrivial smooth family of isometric
embeddings of a closed manifold into Euclidean space with an
optimal dimension. Such a nontrivial family was found for a
polyhedron in Euclidean three space by Connelly
\cite{Connelly1978}. For a general manifold, it is desirable  to
find a canonical isometric embedding into a given Euclidean space
by minimizing the $L^2$ norm of its mean curvature within the
space of isometric embeddings.

Chern told me that he and H. Lewy studied local isometric
embedding of a three manifold into six dimensional Euclidean
space. But they didn't wrote any manuscript on it. The major work
in this subject was done  by E. Berger, Bryant, Griffiths and Yang
\cite{BGY1983} \cite{BBG1983}. They showed that a generic three
dimensional embedding system is strictly hyperbolic, and the
generic four dimensional system is a real principal type. Local
existence is true for a generic metric using a hyperbolic operator
and the Nash-Moser implicit function theorem.

If the target space of isometric embedding is a linear space with
indefinite metric, it is possible that the problem is easier. For
example, by a theorem of Pogorelov \cite{Pogorelov1961-1,
Pogorelov1961-2}, any metric on the two dimensional sphere can be
isometrically embedded into a three dimensional hyperbolic
space-form (where the sectional curvature may be a large negative
constant). Hence it can always be embedded into the hyperboloid of
the Minkowski spacetime. This statement may also be true for
surfaces with higher genus. The fundamental group may cause
obstruction, hence the first step should be an attempt to
canonically embed any complete metric (with bounded curvature) on
a simply connected surface into a three dimensional hyperbolic
space form. It should be also very interesting  to study the
rigidity problem of a space-like surface in Minkowski spacetime.
Besides requesting the metric to be the induced metric, we shall
need one more equation. Such an equation should be related to the
second fundamental form. A candidate appeared in the work of M.
Liu-Yau \cite{LY2003, LY2004} on the quasi-local mass in general
relativity.

In the other direction, Calabi found the condition for a K\"ahler
metric to be isometrically and holomorphically embedded into
Hilbert space with an indefinite signature. In the course of his
investigation, he introduced some kind of distance function that
can be defined by the K\"ahler potential and enjoys many
interesting properties. Calabi's work in this direction which
should be relevant to the flat coordinate appeared in the recent
works of Vafa et al \cite{BCOV1994}.

\begin{quotation}
\noindent {\bf Comment}: {The theory of isometric embedding is a
classical subject. But our knowledge is still rather limited,
especially in dimension greater than three. Many difficult
problems are related to nonlinear mixed type equation or
hyperbolic differential systems over a closed manifold. }
\end{quotation}

\subsection{Rigidity of harmonic maps with negative curvature}

One can define the energy of maps between manifolds and the
critical maps are called harmonic maps. In 1964, Eells-Sampson
\cite{ES1964} and Al'ber \cite{Alb1964} independently  proved the
existence of such maps in their homotopy class if the image
manifold has a non-positive curvature.

When I was working on manifolds with non-positive curvature, I
realized that it is possible to use harmonic map to reprove some
of the theorems in my thesis. I was convinced that it is possible
to use harmonic maps to study rigidity questions in geometry such
as Mostow's theorem \cite{Mos1973}. In 1976, I proved the Calabi
conjecture and applied the newly proved existence of the
K\"{a}hler Einstein metric and the Mostow rigidity theorem to
prove uniqueness of a complex structure on the quotient of the
ball \cite{Yau19773}. Motivated by this theorem, I proposed to use
the harmonic map to prove the rigidity of a complex structure for
K\"{a}hler manifolds with strongly negative sectional curvature. I
proposed this to Siu who carried out the idea when the image
manifold satisfies a stronger negative curvature condition
\cite{Siu1980}. Jost-Yau \cite{JY1983} proved that for harmonic
maps into manifolds with  non-positive curvature, the fibers give
rise to holomorphic foliations even when the map is not
holomorphic. Such a work was found to be used in the work of
Corlette, Simpson et al.

A further result was obtained by Jost-Yau \cite{JY1993} and
Mok-Siu-Yeung \cite{MSY1993} on the proof of the superrigidity
theorem of Margulis \cite{Mar1975}, improving an earlier result of
Corlette \cite{Cor1992} who proved superrigidity for a certain
rank one locally symmetric space. Complete understanding of
superrigidity for the quotient of a complex ball is not yet
available. One needs to find more structures for harmonic maps
which reflect the underlying structure of the manifold. The
analytic proof of super-rigidity was based on an argument of
Matsushima \cite{Mat1962} as was suggested by Calabi. (This was a
topic discussed by Calabi in the special year on geometry in the
Institute for Advanced Study.)

The discrete analogue of harmonic maps is also important. When the
image manifold is a metric space, there are works by Gromov-Schoen
\cite{GS1992}, Korevaar-Schoen \cite{KS1993} and Jost
\cite{Jos1995}. Margulis knew that the super-rigidity for both the
continuous and the discrete case is enough to prove Selberg's
conjecture for the arithmeticity of lattices in groups with rank
$\geq 2$. Unfortunately, the analytic argument mentioned above
only works if the lattices are cocompact as it is difficult to
find a degree one smooth map with finite energy for non-cocompact
lattices. Harmonic maps into a tree have given interesting
applications to group theory. When the domain manifold is a
simplicial complex, there are articles by
Ballmann-\'{S}wwi\c{a}tkowski \cite{BS1997} and M. T. Wang
\cite{Wan1998,Wang20001}, where they introduce maps from complices
which are generalizations of buildings. They also generalized the
work of H. Garland \cite{Gar1973} on the vanishing of the
cohomology group for $p$-adic buildings.

Using the concept of the  center of gravity,
Besson-Courtois-Gallot \cite{BCG1995} give a  metric rigidity
theorem for rank one locally symmetric space. They also proved a
rigidity theorem for manifolds with negative curvature: if the
fundamental group can be split as a nontrivial free product over
some other group $C$, the manifold can be split along a totally
geodesic submanifold with the fundamental group $C$.

\begin{quotation}
\noindent {\bf Comment}: {The harmonic map gives the first step in
matching geometric structures of different manifolds.
Eells-Sampson derived it from the variational principle. One can
also use different elliptic operators to define maps which satisfy
elliptic equations. Higher dimensional applications are mostly
based on the assumption that the image manifold has a metric with
non-positive curvature. In such a case, existence is easier and
uniqueness (as shown by Hartman) is also true. Up to now,
significant results on higher dimensional harmonic maps are based
on such assumptions. Generalization to k\"ahler manifold should be
reasonable. The second  homotopic group should play a role as one
may look at it as a generalization of the work of Sacks-Uhlenbeck.
It may be possible to use harmonic maps to study the  moduli of
geometric structure on a fixed manifold as was done by Michael
Wolf for Riemann surfaces. It will also be nice to see how a
harmonic map can be used to compare graphs.}
\end{quotation}

\subsection{Holomorphic maps}

The works of Liouville, Picard, Schwarz-Pick and Ahlfors show the
importance of hyperbolic complex analysis. Grau\`ert-Reckziegel
\cite{Grauert-Reckziegel1965} generalized this kind of analysis to
higher dimensional complex manifolds. Kobayashi
\cite{Kobayashi1967} and H. Wu \cite{Wu1967} put this theory in an
elegant setting. Kobayashi introduced the concept of hyperbolic
complex manifolds. Its elegant formulation has been influential.
An important application of the negative curvature metric is the
extension theorem for holomorphic maps, as was achieved by the
work of Griffiths-Schmid \cite{Griffiths-Schmid1969} on maps to a
period domain and by  the extension theorem of Borel
\cite{Borel1972} on compactification of Hermitian symmetric space.
A major question was Lang's conjecture: on an algebraic manifold
of a general type, there exists a proper subvariety such that the
image of any holomorphic map from ${\mathbb C}$ must be a subset
of this subvariety. It has deep arithmetic geometric meaning. In
terms of the Kobayashi metric, it says that the Kobayashi metric
is nonzero on a Zariski open set. Many  works were done towards
subvarieties of Abelian variety by Bloch, Green-Griffiths,
Kobayashi-Ochiai, Voitag and Faltings. For generic hypersurfaces
in $\mathbb CP^n$, there is  work by Siu \cite{Siu2002}. They
developed the tool of jet differentials and meromorphic
connections. For algebraic surfaces with $C_1^2>2C_2$, Lu-Yau
\cite{LY1990} proved Lang's conjecture, based on the ideas of
Bogomolov.

\begin{quotation} \noindent {\bf Comment}: {Holomorphic
maps have been studied for a long time. There is no general method
to construct such maps based on the knowledge of topology alone,
except the harmonic map approach proposed by me and   carried out
by Siu, Jost-Yau and others. But the approach is effective only
for manifolds with negative curvature. For rigidity questions, the
most interesting manifolds are K\"{a}hler manifolds with
non-positive Ricci curvature, which give the major chunk of
algebraic manifolds of a general type. The K\"{a}hler-Einstein
metric should provide tools to study such problems. Is there any
intrinsic way, based on the metric, to find the largest subvariety
where the image of all holomorphic maps from the complex line lie?
Deformation theory of such a subvariety should be interesting.
There is also the question of when the holomorphic image of the
complex line will intersect a divisor. Cheng and I did find good
conditions for the complement of a divisor to admit the complete
K\"ahler-Einstein metric. For such a geometry, the holomorphic
line should either intersect the divisor or a subset of some
subvarieties. This kind of questions are very much related to
arithmetic questions if the manifolds are defined over number
fields.}
\end{quotation}

\subsection{Harmonic maps from two dimensional surfaces and pseudoholomorphic curves}

Harmonic maps behave especially well for Riemann surface. Morrey
was the first one who solved the Dirichlet problem for energy
minimizing harmonic map into any Riemannian manifold.

Another major breakthrough was made by Sacks-Uhlenbeck
\cite{SU1981} in 1978 where they constructed minimal spheres in
Riemannian manifolds representing elements in the second homotopy
group using a beautiful extension theorem  of a harmonic map at an
isolated point. By pushing their method further, Siu-Yau
\cite{SiY1980} studied the bubbling process for the harmonic map
and made use of it to prove a stable harmonic map must be
holomorphic under curvature assumptions. As a consequence, they
proved the famous conjecture of Frankel that a K\"{a}hler manifold
with positive bisectional curvature is $\mathbb{C}P^n$, as was
discussed in section [\ref{subsec: 1.3.2}].

Gromov \cite{Gro1985} then realized that a pseudoholomorphic curve
for an almost complex structure can be used in a similar way to
prove rigidity of a symplectic structure on $\mathbb{C}P^n$. The
bubbling process mentioned above was sharpened further to give
compactification of the moduli space of pseudoholomorphic maps by
Ye \cite{Ye1994} and Parker-Wolfson \cite{PW1993}. Based on these
ideas, Kontsevich \cite{Kon1994} introduced the concept of stable
maps and the compactification of their moduli spaces.

The formal definitions of Gromov-Witten invariants and quantum
cohomology were based on these developments and the ideas of
physicists. For example, quantum cohomology was initiated by Vafa
(see, e.g., \cite{Vafa1992}) and his coauthors (the name was
suggested by Greene and me). Associativity in quantum cohomology
was due to four physicists WDVV \cite{Wit1990,DVV1991}. The
mathematical treatment (done by Ruan \cite{Ruan1996} and
subsequently by Ruan-Tian \cite{RT1995}) followed the gluing ideas
of the physicists. Ruan-Tian made use of  the ideas of Taubes
\cite{Tau1984}. But important points were overlooked. A. Zinger
\cite{Zin20021,Zin20022} has recently completed  these arguments.

In close analogy with Donaldson's theory, one needs to introduce
the concept of virtual cycle in the moduli space of stable maps.
The algebraic setting of such a concept is  deeper than the
symplectic case and is more relevant to the development for
algebraic geometry. The major idea was conceived by  Jun Li who
also did the algebraic geometric counterpart of Donaldson's theory
(see \cite{Li1993,LT1998}). (The same comment applies to the
concept of the relative Gromov-Witten invariant, where Jun Li made
the vital contribution in the algebraic setting
\cite{Li20012,Li20022}.) The symplectic version of Li-Tian
\cite{LT19982} ignores difficulties, many of which were completed
recently by A. Zinger \cite{Zin20021,Zin20022}.

Sacks-Uhlenbeck studied harmonic maps from higher genus Riemann
surfaces. Independently, Schoen-Yau \cite{SY1979} studied the
concept of the action of an $L_1^2$ map on the fundamental group
of a manifold. It was used to prove the existence of a harmonic
map with prescribed action on the fundamental group. Jost-Yau
\cite{JY1991} generalized such action on fundamental group to a
more general  setting which allows the domain manifold to be
higher dimensional. Recently F. H. Lin developed this idea further
\cite{Lin19992}. He studied extensively geometric measure theory
on the space of maps (see, e.g., \cite{Lin1989,Lin1999}). The
action on the second homotopy group is much more difficult to
understand. I think there should exist a harmonic map with
nontrivial action on the second homotopic group if such a
continuous map exists. Such an existence theorem will give
interesting applications to K\"{a}hler geometry.

There is a supersymmetric version of harmonic maps studied by
string theorists. This is obtained by coupling the map with Dirac
spinors in different ways (which corresponds to different string
theories). While this kind of world sheet theory is fundamental
for the development of string theory, geometers have not paid much
attention to the supersymmetric harmonic map. Interesting
applications may be found. The most recent paper of Chen, Jost, Li
and Wang \cite{CJLW2004} does address to a related problem where
they studied the regularity and energy identities for
Dirac-Harmonic maps.

\begin{quotation} \noindent {\bf Comment}: {Maps from circle
or Riemann surfaces into a Riemannian manifold give a good deal of
information about the manifold. The capability to construct
holomorphic or pseudo-holomorphic maps from spheres with low
degree was the major reason that Mori, Siu-Yau and Taubes were
able to prove the rigidity of algebraic or symplectic structures
on the complex projective space. It will be desirable to find more
ways to construct such maps from low genus curves to manifolds
that are not of a rational type. Their moduli space can be used to
produce various invariants. An outstanding problem is to
understand the invariants on counting curves of a higher genus
which appeared in the fundamental paper of Vafa et al
\cite{BCOV1994}. }
\end{quotation}

\subsection{Morse theory for maps and topological applications}

The energy functional for maps from $S^2$ into a manifold does not
quite give rise to Morse theory. But the perturbation method  of
Sacks-Uhlenbeck did provide enough information for Micallef-Moore
\cite{MM1988} to prove some structure theorem for manifolds with
positive isotropic curvature. (Micallef and Wang \cite{MW1993}
then proved the vanishing of second Betti number in the even
dimensional case. If the manifold  is irreducible, has
non-negative isotropic curvature and non-vanishing second Betti
number, then they proved that its second Betti number equals to
one and it is K\"{a}hler with positive first Chern class.)

If the image manifold has negative curvature, the theorem of
Eells-Sampson \cite{ES1964} says that any map can be canonically
deformed by the heat flow to a unique harmonic map. Hence the
topology of the space of maps is given by the space of
homomorphism between the fundamental groups of the manifolds. This
gives some information of the topology of manifolds with negative
curvature. Farrell and Jones \cite{FJo1989} have done much deeper
analysis on the differentiable structure of manifolds with
negative curvature.

Schoen-Yau \cite{SY1979} exploited the uniqueness theorem for
harmonic maps to demonstrate that only finite groups can act
smoothly on a manifold which admits a non-zero degree map onto a
compact manifold with negative curvature. The size of the finite
group can also be controlled. If the image manifold has
non-positive curvature, then the only compact continuous group
actions are given by the torus.

The topology of the space of maps into Calabi-Yau manifolds should
be very interesting for string theory. Sullivan
\cite{Sullivan-to-appear} has developed an equivariant homology
theory for loop space. It will be interesting to link such a
theory with quantum cohomology when the manifold has a symplectic
structure.

\begin{quotation} \noindent {\bf Comment}: {Morse theory has
been one of the most powerful tools in geometry and topology as it
connects local to global information. One does not expect full
Morse theory for harmonic maps as we have difficulty even  proving
their existence. However, if their existence can be proven, the
perturbation technique may be used and powerful conclusions may be
drawn. }
\end{quotation}

\subsection{Wave maps}

In early eighties, C. H. Gu \cite{Gu1980} studied harmonic maps
when the domain manifold is the two dimensional Minkowski
spacetime. They are called wave maps. Unfortunately, good global
theory took much longer  to develop as there were not many good a
priori estimates. This subject was studied extensively by
Christodoulou , Klainerman, Tao, Tataru and M. Struwe (see, e.g.,
\cite{CT1993,KS2002,Tao2002,Tat2001,SS1993}). It is hoped that
such theory may shed some light on Einstein equations.

\begin{quotation} \noindent {\bf Comment}: {The geometric or
physical meaning of wave maps should be studied. The problem of
vibrating membrane gives a good motivation to study time-like
minimal hypersurface in a Minkowski spacetime. One can  study the
vibration of submanifold by looking into the minimal time-like
hypersurface with the boundary given by the submanifold. It is a
mystery how such vibrations can be related to the eigenvalues of
the submanifold.}
\end{quotation}

\subsection{Integrable system}

Classically, B\"{a}cklund (1875) was able to find a nonlinear
transformation to create a surface with constant curvature in
$\mathbb{R}^3$ from another one. The nonlinear equation behind it
is the Sine-Gordon equation. Then in 1965, Kruskal and Zabusky
(see \cite{Kru1978}) discovered solitons and subsequently in 1967,
Gardner, Greene, Kruskal and Miura \cite{GGKM1974} discovered the
inverse scatting method to solve the KdV equations. The subject of
a completely integrable system became popular.

Uhlenbeck \cite{Uhl1989} used techniques from integrable systems
to construct harmonic maps from $S^2$ to $U(n)$, Bryant
\cite{Bry1982} and Hitchin \cite{Hit1990} also contributed to
related constructions using twistor theory and spectral curves.
These inspired Burstall, Ferus, Pedit and Pinkall \cite{BFPP1993}
to construct harmonic maps from a torus to any compact symmetric
space. In a series of papers, Terng and Uhlenbeck
\cite{Ter1998,TU2000} used loop group factorizations to solve the
inverse scattering problem and to construct B\"{a}cklund
transformations for soliton equations, including Schr\"{o}dinger
maps from $\mathbb{R}^{1,1}$ to a Hermitian symmetric space. There
are recent attempts by Martin Schmidt \cite{Sch2002} to use an
integrable system to study the Willmore surface.

The integrable system also appeared naturally in several geometric
questions such as the Schottky problem (see Mulase
\cite{Mulase1988}) and the Witten conjecture on Chern numbers of
bundles over moduli space of curves.

Geroch found the Backlund transformation for axially symmetric
stationary solutions of Einstein equations. It will be nice to
find such nonlinear transformations for more general geometric
structures.

\begin{quotation} \noindent {\bf Comment}: {It is always
important to find an explicit solution to nonlinear problem.
Hopefully an integrable system can help us to understand general
structures of geometry.}
\end{quotation}

\subsection{Regularity theory}

The major work on regularity theory of harmonic maps in higher
dimensions was done by Schoen-Uhlenbeck \cite{SU1982,SU1983}.
(There is  a weaker version due to Giaquinta-Giusti \cite{GG1984}
and also the earlier work of Ladyzhenskaya-Ural'ceva and
Hildebrandt-Kaul-Widman where the image manifolds for the maps are
more restrictive.) Leon Simon (see \cite{Sim1996})
 made a deep contribution to the structure
of harmonic maps or minimal subvarieties near their singularity.
This was followed by F. H. Lin \cite{Lin1999}. The following is
still a fundamental problem: Are singularities of harmonic maps or
minimal submanifolds stable when we perturb the metric of the
manifolds? Presumably some of them are. Can we characterize them?
How big is the codimension of generic singularities?

In the other direction Schoen-Yau \cite{SY1978} also proved that
degree one harmonic maps are one to one if the image surface has a
non-positive curvature. Results of this type work only for two
dimensional surfaces. It will be nice to study the set where the
Jacobian vanishes.

\begin{quotation} \noindent {\bf Comment}: {There is a very
rich theory of stable singularity for smooth maps. However, in
most problems, we can only afford to deform certain background
geometric structures, while the extremal objects are still
constrained by the elliptic variational problem. Understanding
this kind of stable singularity should play fundamental roles in
geometry. }
\end{quotation}

\newpage

\section{Submanifolds defined by variational principles}
\label{sec 3}                                               

\subsection{Teichm\"uller space}

The totality of the pair of  polarized K\"{a}hler manifolds with a
homotopic equivalence to a fixed manifold gives rise to the
Teichm\"uller space. For an Algebraic curve, this is the classical
Teichm\"uller space. This space  is important for the construction
of the mapping problem for minimal surfaces of a higher genus.

In fact, given a conformal structure on a Riemann surface
$\Sigma$, a harmonic map from $\Sigma$ to a fixed Riemannian
manifold may minimize energy within a certain homotopy class.
However, it may not be conformal and may not be a minimal surface.
In order to obtain a minimal surface, we need to vary the
conformal structure on $\Sigma$ also. Since the space of conformal
structures on a surface is not compact, one needs to make sure the
minimum can be achieved.

If the map $f$ induces an injection on  the fundamental group of
the domain surface, Schoen-Yau  \cite{SY19795} proved the energy
of the harmonic map is proper on the moduli space of conformal
structure on this surface by making use of a theorem of Linda Keen
\cite{Kee1974}. Based on a theory of topology of the $L_1^2$ map,
they proved the existence of incompressible minimal surfaces. As a
product of this argument, it is possible to find a nice exhaustion
function for the Teichm\"{u}ller  space. Michael Wolf
\cite{Wolf1989} was able to use harmonic maps to give a
compactification of Teichm\"{u}ller space which he proved to be
equivalent to the Thurston compactification. S. Wolpert studied
extensively the behavior of the Weil-Petersson metric (see
Wolpert's survey \cite{Wol2003}). A remarkable theorem of Royden
\cite{Roy1970} says  that the Teichm\"{u}ller metric is the same
as the Kobayashi metric. C. McMullen \cite{McM2000} introduced a
new K\"{a}hler metric on the moduli space which can be used to
demonstrate that the moduli space is hyperbolic in the sense of
Gromov \cite{Gro1991}. The great detail of comparison of various
intrinsic metrics on the Teichm\"uller space had been a major
problem \cite{Yau1986}. It was accomplished recently in the works
of Liu-Sun-Yau \cite{LSY2004,LSY20042}. Actually Liu-Sun-Yau
introduced new metrics with bounded negative curvature and
geometry and found the stability of the logarithmic cotangent
bundle of the moduli spaces. Recently L. Habermann and J. Jost
\cite{HJ2005,HJ20052} also study the geometry of the
Weil-Petersson metric associated to the Bergmann metric on the
Riemann surface instead of the Poincar\'{e} metric.

\begin{quotation} \noindent {\bf Comment}: {For a conformally
invariant variational problem, Teichm\"uller space plays a
fundamental role. It covers the moduli space of curves and in many
ways behaves like a Hermitian symmetric space of noncompact type.
Unfortunately, there is no good canonical realization of it as a
pseudo-convex domain in Euclidean space. For example, we do not
know whether it can be realized as a smooth domain or not.

There is also Teichm\"uller space for other algebraic manifolds,
such as Calabi-Yau manifolds. It is an important question in
understanding their global behavior.}
\end{quotation}

\subsection{Classical minimal surfaces in Euclidean space}

There is  a long and rich history of  minimal surfaces in
Euclidean space. Recent contributions include works by Meeks,
Osserman, Lawson, Gulliver, White, Hildebrandt, Rosenberg, Collin,
Hoffman, Karcher, Ros, Colding, Minicozzi, Rodr\'{i}guez,
Nadirashvili and others (see the reference in Colding and
Minicozzi's survey \cite{CM20052}) on embedded minimal surfaces in
Euclidean space. They come close to classifying complete embedded
minimal surfaces and a good understanding of complete minimal
surface in a bounded domain. For example, Meeks-Rosenberg
\cite{MR2005} proved that the plane and helicoid are the only
properly embedded simply connected minimal surfaces in
$\mathbb{R}^3$.

Calabi also initiated the study of isometric embedding of Riemann
surfaces into $S^N$ as minimal surfaces. The geometry of minimal
spheres and minimal torus was then pursued by many geometers
\cite{Carmo-Wallach1969}, \cite{Chern1970}, \cite{Bry1982},
\cite{Hit1990}, \cite{Lawson1970}.

\begin{quotation} \noindent {\bf Comment}: {This is one of the most beautiful
subjects in geometry where Riemann made important contributions.
Classification of complete minimal surface is nearly accomplished.
However a similar problem for compact minimal surfaces in $S^3$ is
far from being solved. It is also difficult to detect which set of
disjoint Jordan curves can bound a connected minimal surface. The
classification of moduli space of complete minimal surfaces with
finite total curvature should  be studied in detail.}
\end{quotation}

\subsection{Douglas-Morrey solution, embeddedness and application to topology of three manifolds}

In a series of papers started in 1978, Meeks-Yau
\cite{MY1980,MY1981,MY1982,MY19822} settled a classical conjecture
that the Douglas solution for the Plateau problem is embedded if
the boundary curve is a subset of a mean convex boundary. (One
should note that Osserman \cite{Oss1970} had already settled the
old problem of non-existence of branched points for the Douglas
solution while Gulliver \cite{Gul1973} proved non-existence of
false branched points.) We made use of the area minimizing
property of minimal surfaces  to prove these surfaces are
equivariant with respect to the group action. Embedded surfaces
which are equivariant play important roles for finite group
actions on manifolds. Coupling with a theorem of Thurston, we can
then prove the Smith conjecture \cite{YM1984} for cyclic groups
acting on the spheres: that the set of fixed points is not a
knotted curve.

The Douglas-Morrey solution of the Plateau problem is obtained by
fixing the genus of the surfaces. However, it is difficult to
minimize the area when the genus is allowed to be arbitrary large.
This was settled by Hardt-Simon \cite{HS1979} by proving the
boundary regularity of the varifold solution of the Plateau
problem. In the other direction, Almgren-Simon \cite{AS1979}
succeeded in minimizing the area among embedded disks with a given
boundary in Euclidean space. The technique was used by
Meeks-Simon-Yau \cite{MLY1982} to prove the  existence of embedded
minimal spheres enclosing a fake ball. This theorem has been
important to prove that the universal covering of an irreducible
three manifold is irreducible. They also gave conditions for the
existence of embedding minimal surfaces of a higher genus. This
work was followed by topologists Freedman-Hass-Scott
\cite{FHS1983}. Pitts \cite{Pit1981} used the mini-max argument
for varifolds to prove the existence of an embedded minimal
surfaces. Simon-Smith (unpublished) managed to prove the existence
of an embedded minimax sphere for any metric on  the three sphere.
J. Jost \cite{Jost1989} then extended it to find four mini-max
spheres. Pitts-Rubinstein (see, e.g., \cite{PR1988}) continued to
study such mini-max surfaces. Since such mini-max surfaces have
Morse index one, I was interested in representing such minimal
surface as a Heegard splitting of the three manifolds. I estimated
its genus based on the fact that the second eigenvalue of the
stability operator is nonnegative. This argument (dates back to
Szego-Hersch) is to map the surface conformally to $S^2$. Hence we
can use three coordinate functions, orthogonal to the first
eigenfunction, to be trial functions. The estimate gave an upper
bound of the genus for mini-max surfaces in compact manifolds with
a positive scalar curvature. About twenty years ago, I was hoping
to use such an estimate to control a Heegard genus as a way to
prove Poincar\'{e} conjecture. While the program has not
materialized, three manifold topologists did adapt the ideas of
Meeks-Yau to handle combinational type minimal surfaces and gave
applications in three manifold topology.

The most recent works of Colding and Minicozzi
\cite{CM20041,CM20042,CM20043,CM20044}  on lamination by minimal
surfaces and estimates of minimal surfaces without the area bound
are quite remarkable. They \cite{CM2005} made contributions to
Hamilton's Ricci flow by bounding the total time for evolution.
Part of the idea came from the above mentioned inequality.

\begin{quotation} \noindent {\bf Comment}: {The
application of minimal surface theory to three manifold topology
is a very rich subject. However, one needs to have a deep
understanding of the construction of minimal surfaces. For
example, if minimal surfaces are constructed by the method of
mini-max, one needs to know the relation of their Morse index to
the dimension of the family of surfaces that we use to perform the
procedure of mini-max. A detail understanding may lead to a new
proof of the Smale conjecture, as we may construct a minimal
surface by a homotopic group of embeddings of surfaces.
Conversely, topological methods should help us to classify closed
minimal surfaces. }
\end{quotation}

\subsection{Surfaces related to classical relativity}\label{subsec: 3.4}

Besides minimal surfaces, another important class of surfaces are
surfaces with constant mean curvature and also surfaces that
minimize the $L^2$-norm of the mean curvature. It is important to
know the existence of such surfaces in a three dimensional
manifold with nonnegative scalar curvature, as they are relevant
to question in general relativity.

The existence of minimal spheres is related to the existence of
black holes. The most effective method was developed by Schoen-Yau
 \cite{SY1983} where they \cite{SY19793} proved the existence theorem
for the equation of Jang. It should be nice to find new methods to
prove existence of stable minimal spheres. The extremum of the
Hawking mass is related to minimization of the $L^2$ norm of mean
curvature. Their existence and behavior have  not been understood.

For surfaces with constant mean curvature, we have the concept of
stability. (Fixing the volume it encloses, the second variation of
area is non-negative.) Making use of my work on eigenvalues with
Peter Li, I proved with Christodoulou \cite{CY1988} that the
Hawking mass of such a surface is positive. (This was part of my
contribution to the proposed joint project with
Christodoulou-Klainerman which did not materialize.) This fact was
used by Huisken and me \cite{HY1996} to prove uniqueness and the
existence of foliation by constant mean curvature spheres for a
three dimensional asymptotically flat manifold with positive mass.
(We initiated this research in 1986. Ye studied our work and
proved existence of similar foliations under various conditions,
see \cite{Ye1996}.)

This foliation was used by Huisken and Yau \cite{HY1996} to give a
canonical coordinate system at infinity. It defines the concept of
center of gravity where important properties for general
relativity are found. The most notable is that total linear
momentum is equal to the total mass multiple with the velocity of
the center of the mass. One expects to find good asymptotic
properties of the tensors in general relativity along these
canonical surfaces. We hope to find a good definition of angular
momentum based on this concept of center of gravity so that global
inequality like total mass can dominate the square norm of angular
momentum.

The idea of using the foliation of surfaces satisfying various
properties (constant Gauss curvature, for example) to study three
manifolds in general relativity is first developed by R. Bartnik
\cite{Bar1993}. His idea of quasi-spherical foliation gives a good
parametrization of a large class of metrics with positive scalar
curvature.

Some of these ideas were used by Shi-Tam \cite{ST2002} to study
quantities associated to spheres which bound three manifolds with
positive scalar curvature. Such a quantity is realized to be the
quasi-local mass of Brown-York \cite{BY1993}. At the same time,
Melissa Liu and Yau \cite{LY2003,LY2004} were able to define a new
quasi-local mass for general spacetimes in general relativity,
where some of the ideas of Shi-Tam were used. Further works by M.
T. Wang and myself generalized Liu-Yau's work by studying
surfaces in hyperbolic space-form.

My interest in quasi-local mass dates back to the paper that I
wrote with Schoen \cite{SY1983} on the existence of a black hole
due to the condensation of matter. It is desirable to find a
quasi-local mass which includes the effect of matter and the
nonlinear effect of gravity. Hopefully one can prove that when
such a mass is larger than a constant multiple of the square root
of the area, a black hole forms. This has not been achieved.

\begin{quotation} \noindent {\bf Comment}: {When surfaces
theory appears in general relativity, we gain intuitions from both
geometry and physics together. This is a fascinating subject. }
\end{quotation}

\subsection{Higher dimensional minimal subvarieties}

Higher dimensional minimal subvarieties  are very important for
geometry. There are works by Federer-Fleming \cite{FF1960},
Almgren \cite{Alm2002} and Allard \cite{All1975}. The attempt to
prove the Bernstein conjecture, that  minimal graphs are linear,
was a strong drive for its development. Bombieri, De Giorgi and
Giusti \cite{BGG1969} found the famous counterexample to the
Bernstein problem. It initiated a great deal of interest in the
area minimizing cone (as a graph must be area minimizing).
Schoen-Simon-Yau \cite{SSY1975} found a completely different
approach to the proof of Bernstein problem in low dimensions. This
paper on stable minimal hypersurfaces initiated many developments
on curvature estimates for the codimension one stable
hypersurfaces in higher dimension. There are also works by L.
Simon with Caffarelli and Hardt \cite{CHS1984} on constructing
minimal hypersurfaces by deforming stable minimal cones. Recently
N. Wickramasekera \cite{Wic2004,Wic2005}  did some deep work on
stable minimal (branched) hypersurfaces which generalizes
Schoen-Simon-Yau.

Michael-Simon \cite{Michael-Simon1973} proved the Sobolev
inequality and mean valued inequalities for such manifolds. This
enables one to apply the classical argument of harmonic analysis
to minimal submanifolds. For a minimal graph, Bombieri-Giusti
\cite{BG1972} used ideas of De Giorgi-Nash to prove gradient
estimates of the graph. N. Korevaar \cite{Kor1986} was able to
reprove this gradient estimate based on the maximal principle.

The best regularity result for higher codimension was done by F.
Almgren \cite{Alm2002} when he proved that for any area minimizing
variety, the singular set has the  codimension of at least two.
How such a result can be used for geometry remains to be seen.

It was observed by Schoen-Yau \cite{SY19795} that for a closed
stable minimal hypersurface in a manifold with positive scalar
curvature, the first eigenfunction of the second variational
operator can be used to conformally deform the metric so that the
scalar curvature is positive. This provides an induction process
to study manifolds with a positive scalar curvature. For example,
if the manifold admits a nonzero degree map to the torus, one can
then construct stable minimal hypersurfaces inductively until we
find a two dimensional surface with higher genus which cannot
support a metric with positive scalar curvature. At this moment,
the argument encounters difficulty for dimensions greater than
seven as we may have problems of singularity. In any case, we did
apply the argument to prove the positive action conjecture in
general relativity. The question of which type of singularities
for minimal subvariety are generic under metric perturbation
remains a major question for the theory of minimal submanifolds.

Perhaps the most important possible application of the theory of
minimal submanifolds is the Hodge conjecture: whether a multiple
of a $(p,p)$ type integral cohomology class in a projective
manifold can be represented by an algebraic cycle. Lawson made an
attempt by combining a result of Lawson-Simons \cite{LS1973} and
work of J. King \cite{King1971} and Harvey-Shiffman \cite{HS1974}.
(Lawson-Simons proved that currents in $\mathbb{C}P^n$ which are
minimizing with respect to the projective group action are complex
subvarieties.) The problem of how to use the hypothesis of $(p,p)$
type has been difficult. In general, the algebraic cycles are not
effective. This creates difficulties for analytic methods. The
work of King \cite{King1971} and Shiffman \cite{Shi1986} on
complex currents may be relevant.

Perhaps one should generalize the Hodge conjecture to include
general $(p,q)$ classes, as it is possible that every integral
cycle in $\bigoplus_{i=-k}^{k}H^{p-i,p+i}$ is rationally
homologous to an algebraic  sum of minimal varieties such that
there is a $p-k$ dimensional complex space in the tangent space
for almost every point of the variety: it may be important to
assume the metric to be canonical, e.g. the K\"{a}hler Einstein
metric.

A dual question is how to represent a homology class by Lagrangian
cycles which are minimal submanifolds also. When the manifold is
Calabi-Yau, these are special Lagrangian cycles. Since they are
supposed to be dual to holomorphic cycles, there should be an
analogue of the Hodge conjecture. For example if
$\dim_{\mathbb{C}} M=n$ is odd, any integral element in
$\bigoplus_{i+j=n}H^{i,j}$ should be representable by special
Lagrangian cycles up to a rational multiple provided the cup
product of it with the K\"{a}hler class is zero.

A very much related question is:  if the Chern classes of a
complex vector bundle are of $(p,p)$ type, does the vector bundle,
after adding a holomorphic vector bundle, admit a holomorphic
structure? If the above generalization of the Hodge conjecture
holds, there should be a similar generalization for the vector
bundle. It should also be noted that Voisin \cite{Voi2002}
observed that Chern classes of all holomorphic bundles do not
necessarily generate all rational $(p,p)$ classes. On the other
hand, the K\"{a}hler manifold that she constructed is not
projective.

These questions had a lot more success for four dimensional
symplectic manifolds by the work of Taubes both on the existence
of pseudoholomorphic curves \cite{Tau2000} and on the existence of
anti-self-dual connections \cite{Tau1982,Tau1984}. On a K\"{a}hler
surface, anti-self-dual connections are Hermitian connections for
a holomorphic vector bundle. In particular, Taubes gave a method
to construct holomorphic vector bundles over K\"{a}hler surfaces.
Unfortunately this theorem does not provide much information on
the Hodge conjecture as it follows from Lefschetz theorem in this
dimension.

Another important class of minimal varieties is the class  of
special Lagrangian cycles in Calabi-Yau manifolds. Such cycles
were first developed by Harvey-Lawson \cite{HL1982} in connection
to calibrated geometry. Major works were done by Schoen-Wolfson
\cite{SW2001}, Yng-Ing Lee \cite{Lee2003} and Butscher
\cite{But2004}. One expects Lagrangian cycles to be mirror to
holomorphic bundles and special Lagrangian cycles to be mirror to
Hermitian-Yang-Mills connections. Hence by the
Donaldson-Uhlenbeck-Yau theorem, it is related to stability. The
concept of stability for Lagrangian cycles was discussed by Joyce
and Thomas. Since the Yang-Mills flow for Hermitian connection
exists for all time, Thomas-Yau \cite{TY2002} suggested an analogy
with the mean curvature flow for Lagrangian cycles. For stable
Lagrangian cycles, mean curvature flow should converge to special
Lagrangian cycles. See M. T. Wang \cite{Wan2001,Wan20011}, Smoczyk
\cite{Smo2004} and Smoczyk-Wang \cite{SW2002}.  The geometry of
mirror symmetry was explained by Strominger-Yau-Zaslow in
\cite{SYZ1996} using a family of special Lagrangian tori. There
are other manifolds with special holonomy group. They have similar
calibrated submanifolds. Conan Leung has contributed to studies of
such manifolds and their mirrors (see, e.g.
\cite{Leu2002,Leu20022}).

Submanifolds of space forms are called isoparametric if the normal
bundle is flat and the principal curvatures are constants along
parallel normal fields. These were studied by E. Cartan
\cite{Car1939}. Minimal submanifolds, with constant scalar
curvature are believed to be isoparametric surfaces. There is work
done by  Lawson \cite{Law1969}, Chern-de Carmo-Kobayashi
\cite{CDK1970} and  Peng-Terng \cite{PT1983}. Recently there has
been extensive work by Terng and Thorbergsson (see Terng's survey
\cite{Ter1993} and Thorbergsson \cite{Tho2000}). Terng
\cite{Ter1980} related isometric embedded hyperbolic spaces in
Euclidean space to soliton theory. A nice theory of Lax pair and
loop groups related to geometry has been developed.

\begin{quotation} \noindent {\bf Comment}: {The theory of higher
dimensional minimal submanifolds is one of the deepest subjects in
geometry. Unfortunately our knowledge of the subject is not mature
enough to give  applications to solve outstanding problems in
geometry, such as the Hodge conjecture. But the future is bright.
}
\end{quotation}

\subsection{Geometric flows}

The major geometric flows are flows of submanifold driven by mean
curvature, gauss curvature, inverse mean curvature. Flows that
change geometric structures are Ricci flows and Einstein flow.

Mean curvature flow for varifolds was initiated by Brakke
\cite{Bra1978}. The level set approach was studied by many people:
S. Osher,  L. Evans, Giga, etc (see
\cite{Osh1993,ES1991,CGG1991}). Huisken \cite{Hui1984,Hui1986} did
the first important work when the initial surface is convex. His
recent work with Sinestrari \cite{HS19991,HS19992} on mean convex
surfaces is remarkable and gives a good understanding of the
structure of singularities of mean curvature flow. Mean curvature
flow has many geometric applications. For example, the work of
Huisken-Yau mentioned in \ref{subsec: 3.4} was achieved by mean
curvature flow. Mean curvature flow for spacelike hypersurfaces in
Lorentzian manifolds should be very interesting. Ecker
\cite{Eck1997} did interesting work in this direction. It will be
nice to find the Li-Yau type estimate for such flows.

The inverse mean curvature was proposed by Geroch \cite{Ger1973}
to understand the Penrose conjecture relating the mass with the
area of the black hole. Such a procedure was finally carried out
by Huisken-Ilmanen \cite{HI2001} when the scalar curvature is
non-negative. There was a different proof by H. Bray
\cite{Bray2001} subsequently.

Ricci flow has had spectacular successes in recent years. However,
not much progress has been made on the Calabi flow (see Chang's
survey \cite{Cha2005}) for K\"{a}hler metrics. They are higher
order problems where the maximal principle has not been effective.
An important contribution was made by Chr\'{u}sciel \cite{Chr1991}
for Riemann surface. Inspired by the concept of the Bondi mass in
general relativity,  Chru\'sciel was able to give a new estimate
for the Calabi flow. Unfortunately, a higher dimensional analogue
had not been found.

Natural higher order elliptic problems are difficult to handle.
Affine minimal surfaces and Willmore surfaces are such examples.
L. Simon \cite{Sim1986} made an important contribution to the
regularity of the Willmore surfaces. The corresponding flow
problem should be interesting.

The dynamics of Einstein equations for general relativity is a
very difficult subject. The Cauchy problem was considered by many
people: A. Lichnerowicz, Y. Choquet-Bruhat, J. York, V. Moncrief,
H. Friedrich, D. Christodoulou, S. Klainerman, H. Lindblad, M.
Dafermos (see, e.g., \cite{Lic1994}, \cite{CY2001}, \cite{CM2001},
\cite{Fri1996}, \cite{Chr1999}, \cite{KN2004}, \cite{LR2005},
\cite{Daf2003}). But the global behavior is still far from being
understood. The major unsolved problem is to formulate and prove
the fundamental question of Penrose on Cosmic censorship. I
suggested to Klainerman and Christodoulou to consider small
initial data for the Einstein system. The treatment of stability
of Minkowski spacetime was accomplished by
Christodoulou-Klainerman \cite{CK1993} under small perturbation of
flat spacetime and fast fall off conditions. Recently Lindblad and
Rodnianski \cite{LR2005} gave a simpler proof. A few years ago, N.
Zipser (Harvard thesis) added Maxwell equation to gravity and
still proved stability of Minkowski spacetime. There is remarkable
progress on the problem of Cosmic censorship by M. Dafermos
\cite{Daf2003}. He made an important contribution for the
spherical case. Stability for Schwarzschild or Kerr solutions is
far from being known. Finster-Kamran-Smoller-Yau \cite{FKSY2002}
had studied decay properties of Dirac particles with such
background. The work does indicate the stability of these
classical spacetimes.

The no hair theorem for stationary black holes is a major theorem
in general relativity. it was proved by W. Isra\"el
\cite{Israel1967}, B. Carter \cite{Carter1968}, D. Robinson
\cite{Rob1975} and S. Hawking \cite{Hawking1972}. But the proof is
not completely rigorous for the Kerr metric. In any case, the
existing  uniqueness theorem does assume regularity of the horizon
of the black hole. It is not clear to me whether a nontrivial
asymptotically flat solitary solution of a vacuum Einstein
equation has to be the Schwarzschild solution. There is a
possibility that the Killing field is spacelike. In that case,
there may be a new interesting vacuum solution.

There is extensive literature on  spacelike hypersurfaces with
constant mean curvature. The foliation defined by them gives
interesting dynamics of Einstein equation. These surfaces are
interesting even for $\mathbb{R}^{n,1}$. A. Treiberges studied it
extensively \cite{Tre1982}.  Li, Choi-Treibergs \cite{CT1990} and
T. Wan \cite{Wan1992} observed that the Gauss maps of such
surfaces give very nice examples of harmonic maps mapping into the
disk. Recently Fisher and Moncrief used them to study the
evolution equation of Einstein in $2+1$ dimension.

\begin{quotation} \noindent {\bf Comment}: {The dynamics of
submanifolds and geometric structures reveal the true nature of
these geometric objects deeply. In the process of arriving at a
stationary object or a solitary solution, it encounters
singularities. Understanding the structures of such singularity
will solve many outstanding conjectures in topology such as
Shoenflies conjecture.}
\end{quotation}

\newpage

\section{Construction of geometric structures on bundles and manifolds}
\label{sec 4}                                                   

A fundamental question is how to build geometric structures over a
given manifold. In general, the group of topological equivalences
that leaves  this geometric structure invariant should be a
special group. With the exception of symplectic structures, these
groups are usually finite dimensional. When the geometric
structure is unique (up to equivalence), it can be used to produce
key information about the topological structure.

The study of special geometric structures dates back to Sophus
Lie, Klein and Cartan. In most cases, we like to be able to
parallel transport vectors along paths so that we can define the
concept of holonomy group.

\subsection{Geometric structures with noncompact holonomy group}

When the holonomy group is not compact, there are examples of
projective flat structure, affine flat structure and conformally
flat structure. It is not a trivial matter to determine which
topological manifolds admit such structures. Since the structure
is flat, there is a unique continuation property and hence one can
construct a developing map from a suitable cover of the manifold
to the real projective space, the affine space and the sphere
respectively. The map gives rise to a representation of the
fundamental group of the manifold to the real projective group,
the special linear group and the M\"{o}bius group respectively.
This holonomy representation gives a great deal of information for
the geometric structure. Unfortunately, the map is not injective
in general. In the case where it is injective, the manifold can be
obtained as a quotient of a domain by a discrete subgroup of the
corresponding Lie group. In this case, a lot more can be said
about the manifold as the theories of partial differential
equations and discrete groups can play important roles.

\subsubsection{Projective flat structure}

If a projective flat manifold can be projectively embedded as  a
bounded domain, Cheng-Yau \cite{CY1980} were able to construct a
canonical metric from the real Monge-Amp\`{e}re equation which
generalizes the Hilbert metric. When the manifold is two
dimensional, there are works of C. P. Wang \cite{Wan1991} and J.
Loftin \cite{Lof2001} on how to associate such metrics to a
conformal structure and a holomorphic section of the cubic power
of a canonical bundle. This is a beautiful theory related to the
hyperbolic affine sphere mentioned in chapter one.

There are fundamental works by S. Y. Choi, W. Goldman (see the
reference of Choi-Goldman \cite{CG1997}), N. Hitchin
\cite{Hit1992} and others on the geometric decomposition and the
moduli of projective structures on Riemann surfaces. It should be
interesting to extend them to three or four dimensional manifolds.

\subsubsection{Affine flat structure}

It is a difficult question to determine which manifolds admit flat
affine structures. For example, it is still open whether the Euler
number of such spaces is zero, although great progress was made by
D. Sullivan \cite{Sul1976}. W. Goldman \cite{Gol1981} has also
found topological constraints on three manifolds in terms of
fundamental groups. The difficulty arises as  there is no  useful
metric that is compatible with the underlining affine structure.
This motivated Cheng-Yau \cite{CY19801} to define the concept of
affine K\"{a}hler metric.

When Cheng and I considered the concept of affine K\"{a}hler
metric, we thought that it was a natural analogue of K\"{a}hler
metrics. However, compact nonsingular examples are not bountiful.
Strominger-Yau-Zaslow \cite{SYZ1996} proposed the construction of
mirror manifolds by constructing  the  quotient space of a
Calabi-Yau manifold by a special Lagrangian torus. At the limit of
the large K\"{a}hler class, it was pointed out by Hitchin
\cite{Hit1997} that the quotient space admits a natural affine
structure with a compatible affine K\"{a}hler structure. But in
general, we do expect singularities of such structure. It now
becomes a deep  question to understand what kind of singularity is
allowed and how  we build the Calabi-Yau manifold from such
structures. Loftin-Yau-Zaslow \cite{LYZ2004} have initiated the
study of the structure of a "Y" type singularity. Hopefully one
can find an existence theorem for affine structures over compact
manifolds with prescribed singularities along codimension two
stratified submanifolds.

\subsubsection{Conformally flat structure}

Construction of conformally flat manifolds is also a very
interesting topic. Similar to projective flat or affine flat
manifolds, there are simple constraints from curvature
representation for the Pontrjagin classes. The deeper problem is
to  understand the fundamental group and the developing map. When
the structure admits a conformal metric with positive scalar
curvature, Schoen-Yau \cite{SY1988} proved the rather remarkable
theorem that the developing map is injective. Hence such a
manifold must be the quotient of a domain in $S^n$ by a discrete
subgroup of M\"{o}bius transformations. It would be interesting to
classify such manifolds. In this regard, the Yamabe problem as was
solved by Schoen \cite{Sch1984} did provide a conformal metric
with constant scalar curvature. One hopes to be able to use such
metrics to control the conformal structure. Unfortunately the
metric is not unique and a deep understanding of the moduli space
of conformal metrics with constant scalar curvature should be
important.

Kazdan-Warner \cite{KW1975} and Korevaar-Mazzeo-Pacard-Schoen
\cite{KMPS} developed a conformal method to understand Nirenberg's
problem on prescribed scalar curvature. It was followed by
Chen-Lin \cite{CL2001}, Chang-Gursky-Yang \cite{CGY}. Chen-Lin has
related this problem to mean field theory. Their computations in
the relevant degree theory involves deep analysis. One should
generalize their works to functions which are sections of a flat
line bundle because it is related to the previously mentioned work
of Loftin-Yau-Zaslow \cite{LYZ2004}.  In any case, the
integrability condition of Kazdan and Warner is still not fully
understood.

It is curious  that while bundle theory was used extensively in
Riemannian geometry, it has not been used in the study of these
geometries. One can construct real projective space bundles,
affine bundles or sphere bundles by mapping coordinate charts
projectively, affinely or conformally to the corresponding model
spaces (possible with dimensions different from the original
manifold) and gluing the target model spaces together to form
natural bundles. Perhaps one may study their associated
Chern-Simons forms \cite{CS1974}.

Many years ago, H. C. Wang \cite{WangHC1954} proved the theorem
that if a compact complex manifold has trivial holomorphic tangent
bundle, it is covered by a complex Lie group. It will be nice to
generalize and interpret such a theorem in terms of Hermitian
connections on the manifold with a special holonomy group and
torsion.

This program  was discussed in my paper \cite{Yau19932} on
algebraic characterization of locally Hermitian symmetric spaces.
For a holomorphic stable vector bundle $V$, we can form a stable
vector bundles from $V$ by taking irreducible representation of
$GL(n,\mathbb{C})$ from decomposition of  the tensor product
representation $\bigotimes_pV\bigotimes_qV^*$. By twisting with
powers of canonical line bundle, we can form irreducible stable
bundles with trivial determinant line bundle. In general, such
bundles may not have holomorphic sections. If they do, the section
must be parallel with respect to the Hermitian-Yang-Mills
connection on the bundle, and the structure group of $V$ can be
reduced to a smaller group.  Hermitian-Yang-Mills connections with
reduced holonomy group have good geometric properties. We may
formulate a principle: For stable holomorphic bundles, existence
of nontrivial holomorphic  invariants implies the existence of
parallel tensors and therefore the reduction of structure group.
If the holonomy group is reduced to discrete group, the bundle
will provide representations of the fundamental group into unitary
group. This should compare with Wang's theorem when the bundle is
the tangent bundle.

\begin{quotation} \noindent {\bf Comment}: {Geometric
structures with a noncompact holonomy group is less intuitive than
Riemannian geometry. Perhaps we need to deepen our intuitions by
relating them to other geometric structures, especially those
structures that may carry physical meaning. }
\end{quotation}

\subsection{Uniformization for three manifolds}

An important goal of geometry is to build a canonical metric
associated to  a given topology. Besides the uniformization
theorem in two dimensions, the only (spectacular) work in higher
dimensions is the geometrization program of Thurston (see
\cite{Thu1997}).

W. Thurston made use of ideas from Riemann surface theory, W.
Haken's work \cite{Hak1970} on three-manifolds, G. Mostow's
rigidity \cite{Mos1973} to build his manifolds. Many
mathematicians have contributed to the understanding of  this
program  of Thurston's. (e.g., J. Morgan \cite{Morgan1984}, C.
McMullen \cite{McMullen1990, McMullen1992}, J.  Otal
\cite{Otal1996}, J. Porti \cite{Porti1998}.)  Thurston's orbifold
program was finally proved by M. Boileau, B. Leeb and J. Porti
\cite{BLP2001}. However, all this work has to assume that an
incompressible surface (or corresponding surface in case of
orbifold) exists. When R. Hamilton \cite{Ham1982} had his initial
success on his Ricci flow, I suggested (around 1981) to him to use
his flow to break up the manifold and prove Thurston's conjecture.
His generalization of the theory of Li-Yau \cite{LY1986} to Ricci
flow \cite{Ham1988,Ham1993} and his seminal paper in 1996
\cite{Ham1997} on breaking up the manifold mark cornerstone of the
remarkable program. Perelman's recent idea \cite{Per2002,Per2003}
built on these two works and has gone deeply into the problem.
Detailed discussions have been  pursued by Hamilton,
Colding-Minicozzi, Shioya-Yamaguchi, Zhu, Cao and Huisken in the
past two years. Hopefully it may  lead to the final settlement of
the geometrization program.  This theory of Hamilton and Perelman
should be considered as a crowning achievement of geometric
analysis in the past thirty years. Most ideas developed in this
period by geometric analysts are used.

Let me now explain briefly the work of Hamilton and Perelman.

In early 90's, Hamilton \cite{Ham1993,Ham1995,Ham1997} developed
methods and theorems to understand the structure of singularities
of the Ricci flow. Taking up my suggestions, he proved a
fundamental Li-Yau type differential inequality (now called the
Li-Yau-Hamilton estimate) for the Ricci flow with non-negative
curvature in all dimensions. He gave a beautiful interpretation of
the work of Li-Yau and observed the associated inequalities should
be equalities for solitary solutions. He then established a
compactness theorem for (smooth) solutions to the Ricci flow, and
observed (also independently by T. Ivey\cite{Ive1993}) a pinching
estimate for the curvature for three-manifolds. By imposing an
injectivity radius condition, he rescaled the metric to show that
each singularity is asymptotic to one of the three singularity
models. For type I singularities in dimension three, Hamilton
established an isoperimetric ratio estimate to verify the
injectivity radius condition and obtained spherical or neck-like
structures. Based on the Li-Yau-Hamilton estimate, Hamilton showed
that any type II model is either a Ricci soliton with a neck-like
structure or the product of the cigar soliton with the real line.
Similar characterization for type III model was obtained by
Chen-Zhu \cite{CZ2000}. Hence Hamilton had already obtained the
canonical neighborhood structures (consisting of spherical,
neck-like and cap-like regions) for the singularities of
three-dimensional Ricci flow.

But two obstacles remained: one is the injectivity radius
condition and the other is the possibility of forming a
singularity modelled on the product of the cigar soliton with a
real line which could not be removed by surgery. Recently,
Perelman \cite{Per2002} removed these two stumbling blocks in
Hamilton's program by establishing a local injectivity radius
estimate (also called ``Little Loop Lemma" by Hamilton in
\cite{Ham1995}). Perelman proved the Little Loop Lemma in two
ways, one with an entropy functional he introduced in
\cite{Per2002}, the other with a reduced distance function based
on the same idea as Li-Yau's path integral in obtaining their
inequality \cite{Per2002}. This reduced distance question gives
rise to a Gaussian type integral which he called reduced volume.
The reduced volume satisfies monotonicity property. Furthermore,
Perelman \cite{Per2003} developed a refined rescaling argument (by
considering local limits and weak limits in Alexandrov spaces) for
singularities of the Ricci flow on three-manifolds to obtain a
uniform and global version of the canonical neighborhood structure
theorem.

After obtaining the canonical neighborhoods for the singularities,
one performs geometric surgery by cutting off the singularities
and continue the Ricci flow. In \cite{Ham1997}, Hamilton initiated
such a surgery procedure for four-manifolds with a positive
isotropic curvature. Perelman \cite{Per2003} adapted Hamilton's
geometric surgery procedure to three-manifolds. The most important
question is how to prevent the surgery times from accumulations
and make sure there are only a finite number of surgeries on each
finite time interval. When one performs the surgeries with a given
accuracy at each surgery time, it is possible that the errors may
add up which causes the surgery times to accumulate. Hence at each
step of surgery one is required to perform the surgery more
accurate than the former one. In \cite{Per20032}, Perelman
presented a clever idea on how to find ``fine" necks, how to glue
``fine" caps and how to use rescaling arguments to justify the
discreteness of the surgery times. In the process of rescaling for
surgically modified solutions, one encounters the difficulty of
how to use Hamilton's compactness theorem, which works only for
smooth solutions. The idea to overcome such difficulty consists of
two parts. The first part, due to Perelman \cite{Per2003}, is to
choose cutoff radius (in neck-like regions) small enough to push
the surgical regions far away in space. The second part, due to
Chen-Zhu \cite{Chen-Zhu2005} and Cao-Zhu \cite{Cao-Zhu}, is to
show that the solutions are smooth on some uniform small time
intervals (on compact subsets) so that Hamilton's compactness
theorem can be used.

Once surgeries are known to be discrete in time, one can complete
Schoen-Yau's classification \cite{SY1982} for three-manifolds with
positive scalar curvature. For simply connected three manifolds,
if one can show solution to the Ricci flow with surgery extincts
in finite time, Poincar\'{e} conjecture will be proved. Recently,
such a finite extinction time result was proposed by Perelman
\cite{Per20032} and a proof appeared in Colding-Minicozzi
\cite{CM2005}.

For the full geometrization program, one still needs to find the
long time behavior of surgically modified solutions. In
\cite{Hamilton1999}, Hamilton studied the long time behavior of
the Ricci flow on a compact three-manifold for a special class of
(smooth) solutions called ``nonsingular solutions". Hamilton
proved that any (three-dimensional) nonsingular solution either
collapses or subsequently converges to a metric with constant
curvature on the compact manifold, or at large time it admits a
thick-thin decomposition where the thick part consists of a finite
number of hyperbolic pieces and the thin part collapses. Moreover,
by adapting Schoen-Yau's minimal surface arguments to a parabolic
version, Hamilton showed that the boundary of hyperbolic pieces
are incompressible tori. Then by combining with the collapsing
results of Cheeger-Gromov \cite{Cheeger-Gromov1986}, any
nonsingular solution to the Ricci flow is geometrizable.

In \cite{Per2002, Per2003}, Perelman modified Hamilton's arguments
to analyze the long-time behavior of arbitrary solutions to the
Ricci flow and solutions with surgery in dimension three. Perelman
also argued by showing a thick-thin decomposition, except that he
can only show the thin part has (local) lower bound on sectional
curvature. For the thick part, based on Li-Yau-Hamilton estimate,
Perelman established a crucial elliptic type Harnack estimate to
conclude the thick part consisting of hyperbolic pieces. For the
thin part, he announced a new collapsing result which states that
if a three-manifold collapses with a (local) lower bound on the
sectional curvature, it is then a graph manifold. However, the
proof of the new collapsing result has not been published. Shioya
and Yamaguchi \cite{Shioya-Yamaguchi2000,SYa2003} offered a proof
for compact manifolds. Very recently, Cao-Zhu claimed to have a
complete proof for compact manifolds based only on the
Shioya-Yamaguchi's collapsing result.

Hopefully all these arguments can be checked thoroughly in the
near future. It should also be interesting to see whether other
famous problems in three manifold can be settled by analysis: Does
every three dimensional hyperbolic manifold admit a finite cover
with nontrivial first Betti number?

Hyperbolic metrics have been used by topologists to give
invariants for three dimensional manifolds. Thurston
\cite{Thu1982} observed that the volume of a hyperbolic metric is
an important topological invariant. The associated Chern-Simons
\cite{CS1974} invariant, which is defined by mod integers, can be
looked upon as a phase for such manifolds. These invariants
appeared later in Witten's theory of $2+1$ dimensional gravity
\cite{Wit1987} and  S. Gukov \cite{Guk2005} was able to relate
them to fundamental questions in knot theory.

\begin{quotation} \noindent {\bf Comment}: {This is the most
spectacular development in the last thirty years. Once the three
manifold is hyperbolic, Ricci flow does not give much more
information. Happily one may obtain further information by
performing reduction from four dimension Ricci flow to three
dimension by circle action. Is there any effective way to
understand the totality of all hyperbolic manifolds with finite
volume by constructing flows that may break up topology?}
\end{quotation}

\subsection{Four manifolds}

The major accomplishment of Thurston, Hamilton, Perelman et al is
the ability to create a canonical structure on three manifolds.
Such a structure has not even been conjectured for four manifolds
despite the great success of Donaldson invariants and
Seiberg-Witten invariants. Taubes \cite{Tau1992} did prove a
remarkable existence theorem for self-dual metrics on a rather
general class of four dimensional manifolds. Unfortunately their
moduli space is not understood and their topological implication
is not clear at this moment. Since the twistor space of Taubes
metric admits integrable complex structure, ideas from complex
geometry may be helpful. Prior to the construction of Taubes,
Donaldson-Friedman \cite{DF1989} and LeBrun \cite{LeB1991} have
used ideas from twistor theory to construct self-dual metrics on
the  connected sum of $\mathbb{C}P^2$.

The problem of four manifold is the lack of good diffeomorphic
invariants. Donaldson or Seiberg-Witten provide such invariants.
But they are not powerful enough to control the full structure of
the manifold. A true understanding of four manifolds probably
should come from understanding the question of existence of the
integrable complex structures. The Riemann-Roch-Hirzebruch formula
has been the basic tool to find the integrability condition. In
the last twenty-five years, there are nonlinear methods from
K\"{a}her-Einstein metrics, harmonic maps, anti-self-dual
connections and Seiberg-Witten invariants. However, one needs an
existence theorem to find a canonical way to deform an almost
complex structure to an integrable complex structure. What kind of
obstructions do we expect? The work of Donaldson
\cite{Don19962,Don1999} and Gompf \cite{Gom2004} gave a good
characterization of symplectic manifolds in terms of Lefschetz
fibration. It may be useful to know under which condition such
fibration will give rise to complex structures. I did ask several
of my students to work on it. But no definite answer is known. J
Jost and I \cite{JY1992} studied the rigidity part: if a
K\"{a}hler surface has a topological map to a Riemann surface with
higher genus, it can be deformed to be a holomorphic map by
changing the complex structure of the Riemann surface. One can
derive from the work of Griffiths \cite{Gri1971} that every
algebraic surface has a Zariski open set which admits a complete
K\"ahler-Einstein metric with finite volume and is covered by a
contractible pseudo-convex domain. Perhaps one can classify these
manifolds by topological means.

While the Donaldson invariant gave the first counterexample to the
$h$-cobordism theorem and irreducibility (nontrivial connected sum
with manifolds not homotopic to $\mathbb{C}P^2$) of four
manifolds, the Seiberg-Witten invariant gave the remarkable result
that an algebraic surface of general type can not be diffeomorphic
to rational or elliptic surfaces. It also solves the famous Thom
conjecture that holomorphic curves realize the lowest genus for
embedded surface in a K\"{a}hler surface (Kronheimer-Mrowka
\cite{KM1995} and Ozsv\'{a}th-Szab\'{o} \cite{OS2000}). One
wonders whether one can construct a diffeomorphic invariant based
on metrics which are a generalization of K\"{a}hler-Einstein
metrics.

\begin{quotation} \noindent {\bf Comment}: {A good conjectural
statement need to be made on the topology of four manifolds that
may admit an integrable complex structure. Pseudo-holomorphic
curve and fibration by Riemann surfaces should provide important
information. Geometric flows may still be the major tool. }
\end{quotation}

\subsection{Special connections on bundles}

In the seventies, theoretic physicists were very much interested
in the theory of instantons: self-dual connections on four
manifolds. Singer was able to communicate the favor of this
excitement  to the mathematical community which soon led to his
paper with Atiyah and Hitchin \cite{AHS1978} and also the complete
solution of the problem over the four sphere by
Atiyah-Hitchin-Drinfel'd-Manin \cite{AHDM1978} using twistor
technique of Penrose.

 While the paper of Atiyah-Hitchin-Singer \cite{AHS1978} laid the
algebraic and geometric foundation for self-dual connections, the
analytic foundation was laid by Uhlenbeck \cite{Uhl1982,Uhl19822}
where she established the removable singularity theorem and
compactness theorem for Yang-Mills connections. This eventually
led to the fundamental works of Taubes \cite{Tau1984} and
Donaldson \cite{Don1983} which revolutionized four manifold
topology.

In the other direction, Atiyah-Bott \cite{AB1983} applied Morse
theory to the space of connections over Riemann surface. They
solved  important questions on the moduli space of holomorphic
bundles which was studied by Narasimhan, Seshadri, Ramanathan,
Newstead and Harder. In the paper of Atiyah-Bott, Morse theory,
moment map and localization of equivariant cohomology were
introduced on the subject of vector bundle. It laid the foundation
of works in last twenty years.

The analogue of anti-self dual connections over K\"{a}hler
manifolds are Hermitian Yang-Mills connections, which was shown by
Donaldson \cite{Don1985} for K\"{a}hler surfaces and Uhlenbeck-Yau
\cite{UY1986} for general K\"{a}hler manifolds to be equivalent to
the polystability of bundles. (That polystability of bundle is a
consequence of existence of Hermitian Yang-Mills connection was
first observed by   L\"{u}ber \cite{Lub1983}. Donaldson
\cite{Don1987} was able to make use of the theorem of
Mehta-Ramanathan \cite{MRa1984} and ideas of above two papers to
prove the theorem for projective manifold). It was generalized by
C. Simpson \cite{Sim1988}, using ideas of Hitchin \cite{Hit1987},
to bundles with Higgs fields. It has important applications to the
theory of variation of a Hodge structure \cite{Sim19941,Sim1994}.
G. Daskalopoulos and R. Wentworth \cite{DW1997} studied such a
theory for moduli space of vector bundles over curves. Li-Yau
\cite{JLY1987} generalized the existence of Hermitian Yang-Mills
connections to non-K\"{a}hler manifolds. (Buchdahl \cite{Buc1988}
subsequently did the same for complex surfaces.) Li-Yau-Zheng
\cite{LYZ1990} then used the result to give a complete proof of
Bogomolov's theorem for class VII$_{0}$ surfaces. The only missing
parts of the classification of non-K\"{a}hler surfaces are those
complex surfaces with a finite number of holomorphic curves. It is
possible that the argument of Li-Yau-Zheng can be used. One may
want to use Hermitian Yang-Mills connections with poles along such
curves. I expect more applications of Donaldson-Uhlenbeck-Yau
theory to algebraic geometry.

It should be noted that the construction of Taubes \cite{Tau1992}
on the anti-self-dual connection is achieved by singular
perturbation after gluing instantons from $S^4$. The method is
rather different from Donaldson-Uhlenbeck-Yau. While it applies to
arbitrary four manifolds, it does require some careful choice of
Chern classes for the bundle. It will be nice to find a concept of
stability for a  general complex bundle so that a similar
procedure of Donaldson-Uhlenbeck-Yau can be applied. The method of
singular perturbation has an algebraic geometric counterpart as
was found by Gieseker-Li \cite{GLi1994} and O'Grady
\cite{O'Gra1993} who proved that moduli spaces of algebraic
bundles with fixed Chern classes over algebraic surfaces are
irreducible. Li \cite{Lij1995} also obtained information about
Betti number of such moduli space. Not many general theorems are
known for bundles over algebraic manifolds of a higher dimension.
It will be especially useful for bundles over Calabi-Yau
manifolds.

D. Gieseker \cite{Gie1977} developed the geometric invariant
theory for the  moduli space of bundles and introduced the
Gieseker stability of bundles. Conan Leung \cite{Leu1997}
introduced the analytic counterpart of such bundles in his thesis
under my guidance. While it is  a natural concept, there is still
an analytic problem to be resolved. (He assumed the curvature of
the bundles to be uniformly bounded.)

There were attempts by de Bartolomeis-Tian \cite{BT1996} to
generalize Yang-Mills theory to symplectic  manifolds and also by
Tian \cite{Tian2000} to manifolds with a special holonomy group,
as was initiated by the work of Donaldson and Thomas. However, the
arguments for both papers are not complete and still need to be
finished.

For a given natural structure on a manifold, we can often fix a
structure  and linearize the equation to obtain a natural
connection on the tangent bundle. Usually we obtain Yang-Mills
connections with the extra structure given by the holonomy group
of the original structure. It is  interesting to speculate whether
an iterated procedure can be constructed to find an interesting
metric or not. In any case, we can draw analogous properties
between bundle theory and metric theory. The concept of stability
for bundles is reasonably well understood for the holomorphic
category. I believed that for each natural geometric structure,
there should be a concept of stability. Donaldson \cite{Don1987}
was able to explain stability in terms of moment map, generalizing
the work of Atiyah-Bott \cite{AB1983} for bundles over Riemann
surfaces. It will be nice to find  moment maps for other geometric
structures.

\begin{quotation} \noindent {\bf Comment}: {Bundles with
anti-self-dual connections or Hermitian Yang-Mills connections
have been important for geometry. However, we do not have good
estimates of their curvature of such connections. Such an estimate
would be useful to handle important problems such as the
Hartshorne question (see, e.g.,  \cite{BV1974}) on the splitting
of rank two bundle over high dimension complex projective space.}
\end{quotation}

\subsection{Symplectic structures}

Symplectic geometry had many important breakthroughs in the past
twenty years. A moment map was developed by Atiyah-Bott
\cite{AB1984}, Guillemin-Sternberg \cite{GS1982} who proved the
image of the map is a convex polytope. Kirwan and Donaldson had
developed such a theory to be a powerful tool. The
Marsden-Weinstein \cite{MW1974} reduction has become a useful
method in many branches of geometry. At around the same time,
other parts of symplectic topology were developed by Donaldson
\cite{Don1997}, Taubes \cite{Tau1994}, Gompf \cite{Gom1995},
Kronheimer-Mrowka \cite{KM1995} and others.

The phenomenon of symplectic rigidity is manifested by the
existence of symplectic invariants measuring the 2--dimensional
size of a symplectic manifold. The first such invariant was
discovered by Gromov \cite{Gro1985} via pseudo-holomorphic curves.
Hofer \cite{Hof1993} then developed several symplectic invariants
based on variational methods and successfully applied them to
Weinstein conjectures. Ekeland-Hofer \cite{EH1989} introduced a
concept of symplectic capacity and used it to provide a
characterization of a symplectomorphism not involving any
derivatives. The $C^0$-closed property of the symplectomorphism
group as a subgroup of the diffeomorphism group then follows,
which was independently established by Y. Eliashberg
\cite{Eli1987} via wave front methods. Hofer-Zehnder \cite{HZ1994}
introduced another capacity and discovered the displacement-energy
on ${\mathbb R}^{2n}$. By relating the two invariants with the
energy-capacity inequality, Hofer \cite{Hof1990} found a
bi-invariant norm on the infinite dimensional group of Hamiltonian
symplectomorphisms of ${\mathbb R}^{2n}$. The existence of such a
norm has now been established for general symplectic manifolds by
Lalonde-McDuff \cite{LM1995} via pseudo-holomorphic curves and
symplectic embedding techniques. The generalized Weinstein
conjecture on the existence of a periodic orbits of Reeb flows for
many $3$-manifolds including the $3$-sphere was also established
in Hofer \cite{Hof1993} by studying the finite energy
pseudo-holomorphic plane in the symplectization of contact
$3$-manifolds.

Eliashberg-Givental-Hofer \cite{EGH2000} recently introduced the
concept of symplectic field theory, which is about invariants of
punctured pseudo-holomorphic curves in a symplectic manifold with
cylindrical ends. Though it has not been rigorously established,
some applications in contact and symplectic topology have been
found.

By analyzing the singularities of pseudo-holomorphic curves in a
symplectic 4--manifold, D.  McDuff \cite{McD1994} established
rigorously the positivity of intersections of two distinct  curves
and the adjunction formula of an irreducible curve. Applying these
basic properties to symplectic 4-manifolds containing embedded
pseudo-holomorphic spheres with self-intersections at least $-1$,
she was able to construct minimal models of general symplectic
4--manifolds, and classify those containing embedded symplectic
spheres with non-negative self-intersections.

A fundamental question in symplectic geometry is to decide which
topological manifold admits a symplectic structure and how, as was
pointed out by Smith-Thomas-Yau \cite{STY2002}, mirrors of certain
non-K\"{a}hler complex manifolds should be symplectic manifolds.
Based on this point of view, they construct a large class of
symplectic manifolds with trivial first Chern class by reversing
the procedure of Clemens-Friedman on non-K\"{a}hler Calabi-Yau
manifolds \cite{Cle1983, Fri1986}. In dimension four, the Betti
numbers of such manifolds are determined by T. J. Li
\cite{tLi2005}.  In the last ten years, there has been extensive
work  on symplectic manifolds, initiated by Gromov \cite{Gro1985},
Taubes \cite{Tau19961,Tau19962,Tau19991,Tau19992}, Donaldson
\cite{Don19961,Don19962,Don1999} and Gompf \cite{Gom2004}. These
works are based on the understanding of pseudo-holomorphic curves
and Lefschetz fibrations. They are most successful for four
dimensional manifolds. The major tools are Seiberg-Witten theory
\cite{SW19941,SW19942,Wit1994} and analysis. The work of Taubes on
the existence of pseudo-holomorphic curves and the topological
meaning of its counting is one of the deepest works in geometry.
Based on this work, Taubes \cite{Tau19961} was able to prove the
old conjecture that there is only one symplectic structure on the
standard $\mathbb{C}P^2$. However, the following question of mine
is still unanswered:  If $M$ is a symplectic 4-manifold homotopic
to $\mathbb{C}P^2$, is $M$ symplectomorphic to the standard
$\mathbb{C}P^2$?  (The corresponding question for complex geometry
was solved by me in \cite{Yau19773}.) On the other hand, based on
the work of Taubes \cite{Tau1994}, T. J. Li and A. K. Liu
 \cite{LL1995} did  find a wall crossing formula for four dimensional
manifolds that admit metrics with a positive scalar curvature.
Subsequently A. Liu \cite{aLiu1996} gave the classification of
such manifolds. (The surgery result by Stolz \cite{Sto1992} based
on Schoen-Yau-Gromov-Lawson for manifolds with positive scalar
curvature is not effective for four dimensional manifolds.) As
another application of the general wall crossing formula in
\cite{LL1995}, it was proved by T. J. Li and A. Liu in
\cite{LL19952} that there is a unique symplectic structure on
$S^2$-bundles over any Riemann surface. A main result of D. McDuff
in \cite{McD1990} is used here.

McDuff \cite{McD1987} also used a refined bordism type
Gromov-Witten invariant to distinguish  two cohomologous and
deformation equivalent symplectic forms on $S^2\times S^2\times
T^2$, showing that they are not isotopic. Notice that there are
also cohomologous but non-deformation equivalent symplectic forms
on $K3\times S^2$ as shown by Y. Ruan \cite{Rua1994}. In contrast,
it is not known whether examples of this kind exist in dimension
$4$ or not. This phenomenon might be related to the special
features of pseudo-holomorphic curves in a $4$-manifold.

Fukaya and Oh \cite{FO1997} have developed an elaborate theory for
symplectic manifolds with Lagrangian cycles. Pseudo-holomorphic
disks appeared as  trace of motions of  curves according to Floer
theory. Due to boundary bubbles, the Lagrangian Floer homology is
not always defined. Oh \cite{Oh1992-3} developed some works on
pseudo-holomorphic curves with Lagrangian boundary conditions and
extended the Lagrangian Floer homology to all monotone symplectic
manifolds. In order to understand open string theory, Katz-Liu
\cite{KL2001} and Melissa Liu \cite{LiuCC2002} developed the
theory in analogue of the Gromov-Witten invariant for a
holomorphic curve with boundaries on a given Lagrangian
submanifold. Fukaya \cite{Fuk1993} discovered the underlying
$A^{\infty}$ structure of the Lagrangian Floer homology on the
chain level, leading to the Fukaya category. By carefully
analyzing this $A^{\infty}$ structure, Fukaya, Oh, Ohta and Ono in
\cite{FOOO-book} have constructed  a sequence of obstruction
classes which elucidate the rather difficult Lagrangian Floer
homology theory to a great extent. Seidel-Thomas \cite{ST2001} and
W. D. Ruan \cite{Rua2006} discussed Fukaya's category in relation
to Kontsevich's homological mirror conjecture \cite{Kon1995}. One
wonders whether Fukaya's theory can help to construct canonical
metrics for symplectic structures. For symplectic manifolds that
admitan almost complex structure with zero first Chern class, it
would  be nice to construct Hermitian metrics with torsion that
admit parallel spinor. Such structures may be considered as a
mirror to the system constructed by Strominger on non-K\"{a}hler
complex manifolds. Perhaps one can also gain some knowledge by
reduction of $G_2$ or $\text{Spin}(7)$ structures to six
dimensions.

\begin{quotation} \noindent {\bf Comment}: {Geometry from
the  symplectic point of view has seen powerful development in the
past twenty years. Its relation to Seiberg-Witten theory and
mirror geometry is fruitful. More interesting development is
expected.}
\end{quotation}

\subsection{K\"{a}hler structure}

The most interesting geometric structure is the K\"{a}hler
structure. There are two interesting  pre-K\"{a}hler structures.
One is the complex structure and the other is the symplectic
structure. The complex structure is rather rigid for complex two
dimensional manifolds. However it is much more flexible in
dimension greater than two. For example, the twistor space of
anti-self-dual four manifolds admit complex structures. Taubes
\cite{Tau1992} constructed a large class of such manifolds and
hence a large class of complex three manifolds. There is also the
construction of Clemens-Friedman for non-K\"{a}hler Calabi-Yau
manifolds which will be explained later.

For quite a long time, it was believed  that every compact
K\"{a}hler manifold can be deformed to a projective manifold until
C. Voisin \cite{Voi2004,Voi2004-arxiv} found many counterexamples.
We still need to digest the distinction between these two
categories.

Besides some obvious topological obstruction from Hodge theory and
the rational homotopic type theory of
Deligne-Griffiths-Morgan-Sullivan \cite{DGMS1975}, it has been
difficult to decide which complex manifolds admit K\"ahler
structure. Harmonic map argument does give some information. But
it requires the  fundamental group to be large.

Many years ago, Sullivan \cite{Sul19762} proposed to use the
Hahn-Banach theorem to construct K\"ahler metrics. This involves
the  concept of duality and hence closed currents. P. Gauduchon
\cite{Gau1984} has proposed  those Hermitian metrics $\omega$
which is $\partial\bar\partial\omega^{n-1}=0$. Siu \cite{Siu1983}
was able to use these ideas to prove that every $K3$ surface is
K\"ahler. Demailly \cite{Dem1992} did some remarkable work on
regularization of closed positive currents. Singular K\"ahler
metrics have  been studied and used by many researchers. In fact,
in my paper on proving the Calabi conjecture, I proved the
existence of the K\"ahler metrics singular along subvarieties with
control on volume element. They can be used to handle problems in
algebraic geometry, including Chern number inequalities.

\begin{quotation} \noindent {\bf Comment}: {The K\"ahler structure
is one of the richest structures in geometry. Deeper understanding
may require some more generalized structure such as a singular
K\"ahler metric or balanced metrics. }
\end{quotation}

\subsubsection{Calabi-Yau manifolds}

The construction of Calabi-Yau manifolds was based on the
existence of a complex structure which can support a K\"{a}hler
structure and a pluriharmonic volume form.

A fundamental question is whether an almost complex manifold
admits an integrable complex structure when complex dimension is
greater than two. The condition that the first Chern class is zero
is equivalent to the existence of pluriharmonic volume for
K\"{a}hler manifolds. Such a condition is no more true for
non-K\"{a}hler manifolds. It would  be nice to find a topological
method to construct an integrable complex structure with
pluriharmonic volume form.

Once we have an integrable complex structure, we can start to
search for Hermitian metrics with special properties. As was
mentioned earlier, if we would like the geometry to have
supersymmetries, then a K\"{a}hler metric is the only choice if
the connection is torsion free. Further supersymmetry would then
imply the manifold to be Calabi-Yau. However if we do not require
the connection to be torsion free, Strominger \cite{Str1986} did
derive a set of equations  that exhibit supersymmetries without
requiring the manifold to be K\"{a}hler. It is a coupled system of
Hermitian Yang-Mills connections with Hermitian metrics. Twenty
years ago, I tried to develop such a coupled system. The attempt
was unsuccessful as I restricted myself to K\"{a}hler  geometry.
My student Bartnik with Mckinnon \cite{BM1988} did succeed in
doing so in the Lorentzian case. They found non-singular solutions
for such a coupled system. (The mathematically rigorous proof was
provided by Smoller-Wasserman-Yau-Mcleod \cite{SWYM1991} and
\cite{SWY1993}).

The Strominger's system was shown to be solvable in a neighborhood
of a Calabi-Yau structure by Jun Li and myself \cite{JLY2004}. Fu
and I \cite{FY2005} were also able to solve it for certain complex
manifolds which admit no K\"{a}hler structure. These manifolds are
balanced manifolds and were studied by M. Michelsohn
\cite{Mic19822}. These manifolds can be used to  explain some
questions of flux in string theory (see, e.g.,
\cite{BD2002,CCDL}). Since Strominger has shown such manifolds
admits parallel spinors, I have directed my student C.C. Wu to
decompose  cohomology group of such manifolds correspondingly. It
is expected that many theorems in K\"ahler geometry may have
counterparts in such geometry.

Such a structure may help to understand a proposal of Reid
\cite{Rei1987} in connecting Calabi-Yau manifolds with a different
topology. This was initiated by a construction of Clemens
\cite{Cle1983} who proposed to perform complex surgery by blowing
down rational curves with negative normal bundles in a Calabi-Yau
manifold to rational double points. Friedman \cite{Fri1986} found
the condition to smooth out such singularities. Based on this
Clemens-Friedman procedure, one can construct a complex structure
on connected sums of $S^3\times S^3$. It would be nice to
construct Strominger's system on these manifolds.

The Calabi-Yau structure was used by me and others to solve
important  problems in algebraic geometry before it appeared in
string theory. For example,  the proof of the Torelli theorem (by
Piatetskii-Shapiro and Shafarevich \cite{PiS1971}) for a $K3$
surface by Todorov \cite{Tod1980}-Siu \cite{Siu1983} and the
surjectivity of the period map of a $K3$ surface (by Kulikov
\cite{Kul1977}) by Siu \cite{Siu1981}-Todorov \cite{Tod1980} are
important works for algebraic surfaces. The proof of the Bogomolov
\cite{Bog1978}-Tian \cite{Tian1987}-Todorov \cite{Tod1989} theorem
also requires the metric. The last theorem helps us to understand
the moduli space of Calabi-Yau manifolds. It is  important  to
understand the global behavior of the Weil-Petersson geometry for
Calabi-Yau manifolds. C. L. Wang \cite{CLWang1997} was able to
characterize these points which have finite distance in terms of
the metric.

In my talk \cite{Yau1980} in the Congress in 1978, I outlined the
program and the results of classifying noncompact Calabi-Yau
manifolds. Some of this work was written up in Tian-Yau
\cite{TY1990,TY1991} and Bando-Kobayashi
\cite{Bando-Kobayashi1987, Bando-Kobayashi1990}. During the period
of 1984, there was an urgent request by string theorists to
construct Calabi-Yau threefolds with a Euler number equal to $\pm
6$. During the Argonne Lab conference, I \cite{Yau1985}
constructed such a manifold with a $\mathbb{Z}_3$ fundamental
group by taking the quotient of a bi-degree $(1,1)$ hypersurface
in the product of two cubics. Soon afterwards, more examples were
constructed by Tian and myself \cite{TY1987}. However, it was
pointed out by Brian Greene that all the manifolds constructed by
Tian-Yau can be deformed to my original manifold. The idea of
producing Calabi-Yau manifolds by the complete intersection of
hypersurfaces in products of weighted projective space was soon
picked up by Candelas et al \cite{CDLS1986}. By now, on the order
of ten thousand examples of different homotopic types had been
constructed. The idea of using toric geometry for construction was
first performed by S. Roan and myself \cite{RoanYau1987}. A few
years later, the systematic study by Batyrev \cite{Batyrev1994} on
toric geometry allowed one to construct mirror pairs for a large
class of Calabi-Yau manifolds, generalizing the construction of
Greene-Plesser \cite{Greene-Plesser1990} and Candelas et al
\cite{CDLS1986}. Tian and I \cite{TY1987} were also the first one
to apply flop construction to change topology of Calabi-Yau
manifolds. Greene-Morrison-Plesser \cite{GMP1995} then made the
remarkable discovery of isomorphic quantum field theory on two
topological distinct Calabi-Yau manifolds. Most Calabi-Yau
threefolds are a complete intersection of some toric varieties and
they admit a large set of rational curves. It will be important to
understand the reason behind it. Up to now all the Calabi-Yau
manifolds that have a Euler number $\pm 6$ and a nontrivial
fundamental group can be deformed from the birational model of the
manifold (or their mirrors) that I constructed. It would be
important if one could give a proof of this statement.

The most spectacular advancement on Calabi-Yau manifolds come from
the work of Greene-Plesser, Candelas et al on construction of
pairs of mirror manifolds with isomorphic conformal field theories
attached to them. It allows one to calculate Gromov-Witten
invariants. Existence of such mirror pairs  was conjectured by
Lerche-Vafa-Warner \cite{LVW1989} and  rigorous proof of mirror
conjecture was due to Givental \cite{Givental1996} and
Lian-Liu-Yau \cite{LLY1997} independently. The deep meaning of the
symmetry is still being pursued.

In \cite{SYZ1996}, Strominger, Yau and Zaslow proposed a
mathematical explanation for the mirror symmetry conjecture for
Calabi-Yau manifolds. Roughly speaking, mirror Calabi-Yau
manifolds should admit special Lagrangian tori fibrations and the
mirror transformation is a nonlinear analog of the Fourier
transformation along these tori.

This proposal has opened up several new directions in geometric
analysis. The first direction is  the geometry of special
Lagrangian submanifolds in Calabi-Yau manifolds. This includes
constructions of special Lagrangian submanifolds (\cite{Lee2003}
and others) and (special) Lagrangian fibrations by Gross
\cite{Gross1998, Gross1999} and W.D. Ruan
\cite{WDRuan-lagrangian123}, mean curvature of Lagrangian
submanifolds in Calabi-Yau manifolds by Thomas and Yau
\cite{RThomas2001} \cite{TY2002}, structures of singularities on
such submanifolds by Joyce \cite{Joyce2003} and Fourier
transformations along special Lagrangian fibration by
Leung-Yau-Zaslow \cite{LYZ2001} and Leung \cite{Leung2005}.

The second direction is  affine geometry with singularities. As
explained in \cite{SYZ1996}, the mirror transformation at the
large structure limit corresponds to a Legendre transformation of
the base of the special lagrangian fibration which carries a
natural special affine structure with singularities. Solving these
affine problem is not trivial in geometric analysis \cite{Lof2001}
\cite{LYZ2004} and much work is still needed to be done here.

The third direction is  the geometry of special holonomy and
duality and triality transformation in M-theory. In
\cite{GYZ2003}, Gukov, Yau and Zaslow proposed a similar picture
to explain the duality in M-theory. The corresponding differential
geometric structures are fibrations on $G_2$ manifolds by
coassociative submanifolds. These structures are studied by
Kovalev \cite{Kov2003}, Leung and others \cite{Lee-Leung2002}
\cite{LW2004}.

\begin{quotation} \noindent {\bf Comment}: {Although the first
demonstration of the existence of K\"ahler Ricci flat metric was
shown by me in 1976, it was not until the first revolution of
string theory in 1984 that a large group of researchers did
extensive calculations and changed the face of the subject. It is
a subject that provides a good testing ground for analysis,
geometry, physics, algebraic geometry, automorphic forms and
number theory. }
\end{quotation}

\subsubsection{K\"{a}hler metric with harmonic Ricci form and stability}

The existence of a K\"{a}hler Einstein metric with negative scalar
curvature was proved by Aubin \cite{Aub19761} and me
\cite{Yau1978} independently. I \cite{Yau19773} did find its
important applications to solve classical problems in algebraic
geometry, e.g., the uniqueness of complex structure over
$\mathbb{C}P^2$ \cite{Yau19773}, the Chern number inequality of
Miyaoka \cite{Miy1977}-Yau \cite{Yau19773} and the rigidity of
algebraic manifolds biholomorphic to Shimura varieties. The
problem of existence of K\"{a}hler Einstein metrics with positive
scalar curvature in the general case is not solved. However, my
proof of the Calabi conjecture already provided all the necessary
estimates except some integral estimate on the unknown. This of
course can be turned into hypothesis. I conjectured that an
integral estimate of this sort is related to the stability of
manifolds. Tian \cite{Tian1997} called it K-stability. Mabuchi's
functional \cite{Mab1986}  made the integral estimate  more
intrinsic and it  gave rise to a natural variational formulation
of the problem. Siu has pointed out that the work of Tian
\cite{Tian19902} on two dimensional surfaces is not complete. The
work of Nadel \cite{Nad1990} on the multiplier ideal sheaf did
give useful methods for the subject of the  K\"ahler-Einstein
metric.

For K\"{a}hler Einstein manifolds with positive scalar curvature,
it is possible that they admit a continuous group of
automorphisms. Matsushima \cite{Mat1957} was the first one to
observe that such a group must be reductive. Futaki \cite{Fut1983}
introduced a remarkable invariant for general K\"ahler manifolds
and proved that it must vanish for such manifolds. In my seminars
in the eighties, I proposed that Futaki's theorem should be
generalized to understand the projective group acting on the
embedding of the manifold by a high power of anti-canonical
embedding and that Futaki's invariant should be relevant to my
conjecture \cite{Yau1993} relating the K\"{a}hler Einstein
manifold to stability. Tian asked what happens when manifolds have
no group actions. I explained that the shadow of the group action
is there once it is inside the projective space and one should
deform the manifold to a possibly singular variety to obtain more
information. The connection of Futaki invariant to stability of
manifolds has finally appeared in the recent work of Donaldson
\cite{Don2001,Don2002}. Donaldson introduced a remarkable concept
of stability based purely on concept of algebraic geometry. It is
not clear that Donaldson's algebraic definition has anything to do
with Tian's analytic definition of stability. Donaldson proved
that the existence of K\"{a}hler-Einstein did imply his
K-stability which in turn implies Hilbert stability and asymptotic
Chow stability of the manifold. This theorem of Donaldson already
gives nontrivial information for manifolds with negative first
Chern class and Calabi-Yau manifolds, where existence of
K\"{a}hler-Einstein metrics was proved. Some part of the deep work
of Gieseker \cite{Gie1977} and Viehweg \cite{Viehweg1995} can be
recovered by these theorems. One should also mention the recent
interesting work of Ross-Thomas \cite{RT20041,RT20042} on the
stability of manifolds. Phong-Strum \cite{PS2004} also studied
solutions of certain degenerate Monge-Amp\'{e}re equations and
\cite{PS2005} the convergence of the K\"{a}hler-Ricci flow.

A K\"{a}hler metric with constant scalar curvature is equivalent
to the harmonicity of the first Chern form. The uniqueness theorem
for harmonic K\"{a}hler metric was due to X. Chen \cite{Che2000},
Donaldson \cite{Don2001} and Mabuchi for various cases. (Note that
the most important case of the uniqueness of the K\"{a}hler
Einstein metric with positive scalar curvature was due to the
remarkable argument of Bando-Mabuchi \cite{BM1987}.) My general
conjecture for existence of harmonic K\"{a}hler manifolds based on
stability of such manifolds is still largely unknown. In my
seminar in the mid-eighties, this problem was discussed
extensively. Several students of mine, including Tian
\cite{Tian1990}, Luo \cite{Luo1998} and Wang \cite{xWan2002} had
written thesis related to this topic. Prior to them, my former
students Bando \cite{Ban1984} and Cao \cite{Cao1985} had made
attempts to study the problem of constructing K\"{a}hler-Einstein
metrics by Ricci flow. The fundamental curvature estimate was due
to Cao \cite{Cao1992}. The K\"{a}hler Ricci flow may either
converge to K\"{a}hler Einstein metric or K\"{a}hler solitons.
Hence in order for the approach, based on Ricci flow, to be
successful, stability of the projective manifold should be related
to  such K\"{a}hler solitons. The study of harmonic K\"{a}hler
metrics with constant scalar curvature on toric variety was
initiated by S. Donaldson \cite{Don2002}, who proposed to study
the  existence problem via the real Monge-Amp\`{e}re equation.
This problem of Donaldson in the K\"ahler-Einstein case was solved
by Wang-Zhu \cite{WZ2004}. LeBrun and his coauthors \cite{KLP1997}
also have found special constructions, based on twistor theory,
for harmonic K\"{a}hler surfaces. Bando was also interesting in
K\"{a}hler manifolds with harmonic $i$-th Chern form. (There
should be an analogue of stability of algebraic manifolds
associated to manifolds with harmonic $i$-th Chern form.)

In the early 90's, S.W. Zhang \cite{Zh1996} studied heights of
manifolds. By comparing metrics on Deligne pairings, he found that
a projective variety is Chow semistable if and only if it can be
mapped by an element of a special linear group to a balanced
subvariety. (Note that a subvariety in $\mathbb{C}\text{P}^N$ is
called balanced if the integral of the moment map with respect to
$SU(N+1)$ is zero, where the measure for the integral is induced
from the Fubini-Study metric.) Zhang communicated his results to
me. It is clearly related to K\"{a}hler-Einstein metric and I
urged my students, including Tian, to study this connection.

Zhang's work has a nontrivial consequence on the previous
mentioned development of Donaldson \cite{Don2001, Don2002}. Assume
the projective manifold is embedded by an ample line bundle $L$
into projective space. If the manifold has a finite automorphism
group and admits a harmonic K\"{a}hler metric in $c_1(L)$, then
Donaldson showed that for $k$ large, $L^k$ gives rise to an
embedding which is balanced. Furthermore, the induced Fubini-Study
form divided by $k$ will converge to the harmonic K\"{a}hler form.
Combining the work of Zhang and Luo, he then proved that the
manifold is given by the embedding of $L^k$ is stable in the sense
of geometric invariant theory. Recently, Mabuchi generalized
Donaldson's theorem to certain case which allow nontrivial
projective automorphism.

Donaldson considered the problem from the point of view of
symplectic geometry (K\"{a}hler form is a natural symplectic
form). The Hamiltonian group then acts on the Hilbert space $H$ of
square integrable sections of the line bundle $L$ where the first
Chern class is the K\"{a}hler form. For each integrable complex
structure on the manifold compatible with the symplectic form, the
finite dimensional space of holomorphic sections gives a subspace
of $H$. The Hamiltonian group acts on the Grassmannian of such
subspaces. The moment map can be computed to be related to the
Bergman kernel $\sum_\alpha s_\alpha(x)\otimes s_\alpha^*(y)$
where $s_\alpha$ form an orthonormal basis of the holomorphic
sections. On the other hand, Fujiki \cite{Fujiki1990} and
Donaldson \cite{Don1997} computed the moment map for the
Hamiltonian group action on the space of integrable complex
structure, which turns out to be the scalar curvature of the
K\"{a}hler metric. These two moment maps may not match, but for
the line bundle $L^k$ with large $k$, one can show that they
converge to each other after normalization. Lu \cite{Lu2000} has
shown the first term of the expansion (in terms of $1/k$) of the
Bergman kernel gives rise to scalar curvature. Hence we see the
relevance of constant scalar curvature for a K\"{a}hler metric to
these witha constant Bergman kernel function. S.W. Zhang's result
says that the manifold is Chow semistable if and only if it is
balanced. The balanced condition implies that there is a
K\"{a}hler metric where the Bergman kernel is constant. With the
work of Zhang and Donaldson, what remains to settle my conjecture
is the convergence of the balanced metric when $k$ is large. In
general, we should not expect this to be true. However, for toric
manifolds, this might be the case.

It may be noted that in my paper with Bourguignon and P. Li
\cite{BLY1994} on giving an upper estimate of the first eigenvalue
of an algebraic manifold, this balanced condition also appeared.
Perhaps first eigenfunction may play a role for questions of
stability.

\begin{quotation} \noindent {\bf Comment}: {K\"ahler metrics
with constant scalar curvature is a beautiful subject as it is
related to structure questions of algebraic varieties including
the concept of stability of manifolds. The most effective
application of such metrics to algebraic geometry are still
restricted to the K\"ahler-Einstein metric. The singular
K\"ahler-Einstein metric as was initiated by my paper on Calabi
conjecture should be studied further in application to algebraic
geometry. }
\end{quotation}

\subsection{Manifolds with special holonomy group}

Besides K\"{a}hler manifolds, there are manifolds with special
holonomy groups. Holonomy groups of Riemannian manifolds were
classified by Berger \cite{Ber1955}. The most important ones are
$O(n)$, $U(n)$, $SU(n)$, $G_2$ and $Spin(7)$. The first two groups
correspond to Riemannian and K\"{a}hler geometry respectively.
$SU(n)$ corresponds to Calabi-Yau manifolds. A $G_2$ manifold is
seven dimensional and a $Spin(7)$ is eight dimensional (assuming
they are irreducible manifolds). These last three classes of
manifolds have zero Ricci curvature. It may be noted that before I
\cite{Yau1978} proved the Calabi conjecture in 1976, there was no
known nontrivial compact Ricci flat manifold. Manifolds with a
special holonomy group admit nontrivial parallel spinors and they
correspond to supersymmetries in the language of physics. The
input of ideas from string theory did give a lot of help to
understand these manifolds. However, the very basic question of
constructing these structures on a given topological space is
still not well understood. In the case of $G_2$ and $Spin(7)$, it
was initiated by Bryant (see \cite{Bry1987,BS1989}). The first set
of compact examples was given by Joyce
\cite{Joy1996,Joy19962,Joy1999}. Recently Dai-Wang-Wei
\cite{DWW2005} proved the stability of manifolds with parallel
spinors.

The nice construction of Joyce was based on a singular
perturbation which is similar to the construction of Taubes
\cite{Tau1982} on anti-self-dual connections. However, it is not
global enough to give a good parametrization of $G_2$ or $Spin(7)$
structures. A great deal more work is needed. The beautiful theory
of Hitchin
 \cite{Hit2000,Hit2001} on three forms and four forms may lead to a
resolution of these important problems.

\begin{quotation} \noindent {\bf Comment}: {Recent interest
in M-theory has stimulated a lot of activities on manifolds with
special holonomy group. We hope a complete structure theorem for
such manifolds can be found. }
\end{quotation}

\subsection{Geometric structures by reduction}
One can also obtain new geometric structures by imposing some
singular  structures on a manifold with a special holonomy group.
For example, if we require a metric cone to admit a $G_2$,
$Spin(7)$ or Calabi-Yau structure, the link of the cone will be a
compact manifold with special structures. They give interesting
Einstein metrics. When the cone is Calabi-Yau, the structure on
the odd dimensional manifold is called Sasakian Einstein metric.

There is a natural Killing field called the Reeb vector field
defined on a Sasakian Einstein manifold. If it generates a circle
action, the orbit space gives rise to a K\"{a}hler Einstein
manifold with positive scalar curvature. However, it need not
generate a circle action and J. Sparks, Gauntlelt, Martelli and
Waldram \cite{GMSW} gave many interesting explicit examples of
non-regular Sasakian Einstein structures. They have interesting
properties related to conformal field theory. For quasi-regular
examples, there was work by  Boyer, Galicki and Koll\'{a}r
\cite{BKK2003}. The procedure  gave many interesting examples of
Einstein metrics on odd dimensional manifolds.

Sparks, Matelli and I have been pursuing general theory of
Sasakian Einstein manifolds. I would like to consider them as a
natural generalization of K\"ahler manifolds.

\begin{quotation} \noindent {\bf Comment}: {The recent
development of the Sasakian Einstein metric show that it gives a
natural generalization of  the K\"ahler-Einstein metric. Its
relation with the recent activities on ADS/CFT theory is
exciting.}
\end{quotation}

\subsection{Obstruction for existence of Einstein metrics on
general manifolds}

The existence of Einstein metrics on a fixed topological manifold
is clearly one of the most important questions in geometry. Any
metrics with a compact special holonomy group are Einstein.
Besides K\"{a}hler geometry, we do not know much of their moduli
space. For an Einstein metric with no special structures, we know
only some topological constraints on four dimensional manifolds.
There is work by Berger \cite{Ber1961}, Gray \cite{Gra1972} and
Hitchin \cite{Hit19742} in terms of inequalities linking a Euler
number and the signature of the manifold. (This is of course based
on Chern's work \cite{Chern1946} on the representation of
characteristic classes by curvature forms.) Gromov \cite{Gro1982}
made use of his concept of Gromov volume to give further
constraint. LeBrun \cite{LeB2001} then introduced the ideas from
Seiberg-Witten invariants to enlarge such classes and gave
beautiful rigidity theorems on Einstein four manifolds.
Unfortunately it is very difficult to understand moduli space of
Einstein metrics when they admit no special structures. For
example, it is still an open question of whether there is only one
Einstein metric on the four dimensional sphere. M. Wang and Ziller
\cite{WZ1986} and C. Boehm \cite{Boehm1998} did use symmetric
reductions to give many examples of Einstein metrics for higher
dimensional manifolds. There may be much more examples of Einstein
manifolds with negative Ricci curvature than we expected. This is
certainly true for compact manifolds, with negative Ricci
curvature. Gao-Yau \cite{Gao-Yau1986} was the first one to
demonstrate that such a metric exists on the three sphere. A few
years later, Lokhamp \cite{Lokhamp1995} used the $h$-principle of
Gromov to prove such a metric exists on any manifold with a
dimension greater than three. It would be nice to prove that every
manifold with a dimension greater than $4$ admits an Einstein
metric with negative Ricci curvature.

\begin{quotation} \noindent {\bf Comment}: {The  Einstein manifold
without extra special structures is a difficult subject. Do we
expect a general classification for such an important geometric
structure?}
\end{quotation}

\subsection{Metric Cobordism}

In the last five years, a great deal of attention was addressed by
physicists on the holographic principle: boundary geometry should
determine the geometry in the interior. The ADS/CFT correspondence
studies the  conformal boundary of the Einstein manifold which is
asymptotically hyperbolic. Gauge theory on the boundary is
supposed to be dual to the theory of gravity in the bulk. Much
fascinating work was done in this direction. Manifolds with
positive scalar curvature appeared as conformal boundary are
important for physics.  Graham-Lee \cite{GLe1991} have studied a
perturbation problem near the standard sphere which bounds the
hyperbolic manifold. Witten-Yau \cite{WY2000} proved that for a
manifold with positive scalar curvature to be a conformal
boundary, it must be connected. It is not known whether there are
further obstructions.

Cobordism theory had been a powerful tool to classify the topology
of manifolds. The first fundamental work was done by Thom who
determined the cobordism group. Characteristic numbers play
important roles. When two manifolds are cobordant to each other,
the theory of surgery helps us to deform one manifold to another.
It is clear that any construction of surgery that may preserve
geometric structures would play a fundamental role in the future
of geometry.

There are many geometric structures that are preserved under a
connected sum construction. This includes the category of
conformally flat structures, metrics with positive isotropic
curvature and metrics with positive scalar curvature. For the last
category, there was work by Schoen-Yau-Gromov-Lawson where they
perform surgery on spheres with a codimension greater than or
equal to $3$. A key part of the work of Hamilton-Perelman is to
find a canonical neighborhood to perform surgery. If we can deform
the spheres in the above SYGL construction to a more canonically
defined position, one may be able able to create an extra
geometric structure for the result of SYGL. In fact, the
construction of Schoen-Yau did provide some information about the
conformal structures of the manifold. In complex geometry, there
are two important canonical neighborhoods given by the log
transform of Kodaira and the operation of flop. There should be
similar constructions for other geometric structures.

The theory of quasi-local mass mentioned in section \ref{subsec:
3.4} is another example of how boundary geometry can be controlled
by the geometry in the bulk. The work of Choi-Wang \cite{CW1983}
on the first eigenvalue is also based on the manifold that it
bounds. There can be interesting theory of metric cobordism.

In the other direction, there are also beautiful rigidity of
inverse problems for metric geometry by Gerver-Nadirashvili
\cite{GN1984} and  Pestov-Uhlmann \cite{PU2005} on recovering  a
Riemannian metric when one knows the distance functions between
pair of points on the boundary, if the Riemannian manifold is
reasonably convex.

\begin{quotation} \noindent {\bf Comment}: {There should be
a mathematical foundation of the holographic principle of
physicists. Good understanding of metric cobordism may be useful.}
\end{quotation}

\end{document}